\documentclass[10pt,twoside]{article}
\def\YEAR{\year}\newcount\VOL\VOL=\YEAR\advance\VOL by-1995
\def\firstpage{1}\def\lastpage{1000}
\def\received{}\def\revised{}
\def\communicated{}

\makeatletter
\def\magnification{\afterassignment\m@g\count@}
\def\m@g{\mag=\count@\hsize6.5truein\vsize8.9truein\dimen\footins8truein}
\makeatother

\font\eightrm=cmr8
\font\caps=cmcsc10          
\font\Caps=cmcsc10 scaled \magstep1  

\makeatletter
\@specialpagefalse\headheight=8.5pt
\def\DocMath{}
\renewcommand{\@evenhead}{%
  \ifnum\thepage>\lastpage\rlap{\thepage}\hfill%
  \else\rlap{\thepage}\slshape\leftmark\hfill{\caps\SAuthor}\hfill\fi}%
\renewcommand{\@oddhead}{%
  \ifnum\thepage=\firstpage{\DocMath\hfill\llap{\thepage}}%
  \else{\slshape\rightmark}\hfill{\caps\STitle}\hfill\llap{\thepage}\fi}%
\makeatother

\def\TSkip{\bigskip}
\newbox\TheTitle{\obeylines\gdef\GetTitle #1
\ShortTitle #2
\SubTitle  #3
\Author   #4
\ShortAuthor #5
\EndTitle
{\setbox\TheTitle=\vbox{\baselineskip=20pt\let\par=\cr\obeylines%
\halign{\centerline{\Caps##}\cr\noalign{\medskip}\cr#1\cr}}%
	\copy\TheTitle\TSkip\TSkip%
\def\next{#2}\ifx\next\empty\gdef\STitle{#1}\else\gdef\STitle{#2}\fi%
\def\next{#3}\ifx\next\empty%
  \else\setbox\TheTitle=\vbox{\baselineskip=20pt\let\par=\cr\obeylines%
  \halign{\centerline{\caps##} #3\cr}}\copy\TheTitle\TSkip\TSkip\fi%
\centerline{\caps #4}\TSkip\TSkip%
\def\next{#5}\ifx\next\empty\gdef\SAuthor{#4}\else\gdef\SAuthor{#5}\fi%
\ifx\received\empty\relax
  \else\centerline{\eightrm Received: \received}\fi%
\ifx\revised\empty\TSkip%
  \else\centerline{\eightrm Revised: \revised}\TSkip\fi%
\ifx\communicated\empty\relax
  \else\centerline{\eightrm Communicated by \communicated}\fi\TSkip\TSkip%
\catcode'015=5}}\def\Title{\obeylines\GetTitle}
\def\Abstract{\begingroup\narrower
  \parskip=\medskipamount\parindent=0pt{\caps Abstract. }}
\def\EndAbstract{\par\endgroup\TSkip}

\long\def\MSC#1\EndMSC{\def\arg{#1}\ifx\arg\empty\relax\else
   {\par\narrower\noindent%
   2000 Mathematics Subject Classification: #1\par}\fi}

\long\def\KEY#1\EndKEY{\def\arg{#1}\ifx\arg\empty\relax\else
	{\par\narrower\noindent Keywords and Phrases: #1\par}\fi\TSkip}

\newbox\TheAdd\def\Addresses{\vfill\copy\TheAdd\vfill
  \ifodd\number\lastpage\vfill\eject\phantom{.}\vfill\eject\fi}
{\obeylines\gdef\GetAddress #1
\Address #2 
\Address #3
\Address #4
\EndAddress
{\def\xs{4.3truecm}\parindent=0pt
\setbox0=\vtop{{\obeylines\hsize=\xs#1\par}}\def\next{#2}
\ifx\next\empty 
   \setbox\TheAdd=\hbox to\hsize{\hfill\copy0\hfill}
\else\setbox1=\vtop{{\obeylines\hsize=\xs#2\par}}\def\next{#3}
\ifx\next\empty 
   \setbox\TheAdd=\hbox to\hsize{\hfill\copy0\hfill\copy1\hfill}
\else\setbox2=\vtop{{\obeylines\hsize=\xs#3\par}}\def\next{#4}
\ifx\next\empty\ 
   \setbox\TheAdd=\vtop{\hbox to\hsize{\hfill\copy0\hfill\copy1\hfill}
        \vskip20pt\hbox to\hsize{\hfill\copy2\hfill}}
\else\setbox3=\vtop{{\obeylines\hsize=\xs#4\par}}
   \setbox\TheAdd=\vtop{\hbox to\hsize{\hfill\copy0\hfill\copy1\hfill}
	    \vskip20pt\hbox to\hsize{\hfill\copy2\hfill\copy3\hfill}}
\fi\fi\fi\catcode'015=5}}\gdef\Address{\obeylines\GetAddress}

\usepackage{amsfonts,amsmath,amssymb,amsopn}

\newtheorem{thm}{Theorem}[section]
\newtheorem{prop}[thm] {Proposition}
\newtheorem{cor}[thm] {Corollary}
\newtheorem{lem}[thm] {Lemma}

\begin{document}

\Title
The overconvergent site I. Coefficients
\ShortTitle 
\SubTitle  
\Author 
Bernard Le Stum
\ShortAuthor 
\EndTitle
\Abstract 
We define and study the overconvergent site of an algebraic variety, the sheaf of overconvergent functions on this site and show that the modules of finite presentations correspond to Berthelot's overconvergent isocrystals.
We work with Berkovich theory instead of rigid analytic geometry and do not use any of Berthelot's results.
This gives a complete alternative approach to rigid cohomology.
\EndAbstract
\MSC 
14F30
\EndMSC
\KEY 
\EndKEY
\Address 
IRMAR
Universit\'e de Rennes 1
Campus de Beaulieu
35042 Rennes Cedex
France

\Address
\Address
\Address
\EndAddress

\tableofcontents

\section*{Introduction}

In order to give an algebraic description of Betti cohomology, one can use de Rham cohomology which can then be interpreted as the cohomology of the infinitesimal site.
The category of coefficients, locally trivial families of finite dimensional vector spaces, is replaced successively by coherent modules with integrable connexions, and then, modules of finite presentation.
In the positive characteristic situation, de Rham cohomology has to be replaced with rigid cohomology and modules with integrable connexions with overconvergent isocrystals.
We will define here the overconvergent site which plays the role of the infinitesimal site.
A first hint at this approach is already in Berthelot's fundamental article (\cite {Berthelot96*}, 2.3.2. ii)) and this is actually the way I liked to define isocrystal in my PhD Thesis (see also \cite {EtesseLeStum93}, 1.1).

Beside its intrinsic interest, there are many reasons to look for such an interpretation of rigid cohomology.
For example, we should get for free a Leray spectral sequence giving the overconvergence of the Gauss-Manin connexion. Also, our setting should be well-suited to describe Besser's integration (\cite {Besser00}) or Chiarellotto-Tsuzuki's descent theory (\cite {ChiarellottoTsuzuki03}).
Actually, using \'etale topology should give an interpretation of the slopes by means of $p$-adic cohomology.
Finally, replacing schemes by log-schemes should give a comparison theorem with log-crystalline cohomology.
In order to avoid technical complications in this first attempt, we should not consider \'etale cohomology nor log-schemes.

In our presentation, we will systematically replace rigid geometry with analytic geometry in the sense of Berkovich.
The main reason is that, in the construction of rigid geometry, strict neighborhoods play an essential role.
In Berkovich theory those are just usual neighborhoods.
Note also that the notion of generic point that is central in Dwork's theory has a very natural interpretation using Berkovich theory (\cite {Chiarellotto00}). Anyway, the reader can always read ``rigid variety'' when we write ``analytic space'', ``wide open subset'' when we write ``open subset'' and ``open subset'' when we write ``analytic domain''.

Of course, this article owes much to Berthelot's previous work on rigid cohomology. We only want to rewrite his theory with a slightly different approach. The reader should note however that we do not make any use of Berthelot's results and that, for this purpose, the article is totally self-contained. In particular, the reader need not know anything about rigid cohomology.

Let us now present our main result in a more precise form.
Let $K$ be a complete ultrametric field of characteristic zero with valuation ring $\mathcal V$ and residue field $k$.
To each pair composed of a formal $\mathcal V$-scheme $S$ and an $S_k$-scheme $X$, we associate in a natural way a site $\textrm{AN}^\dagger(X/S)$.
On this site lives the sheaf of ring $\mathcal O^\dagger_{X/S}$ of dagger functions.
Our main theorem says that the category of $\mathcal O^\dagger_{X/S}$-modules of finite presentation is equivalent to the category of overconvergent isocrystals on $X/S$.

At this point, I should recall the definition of the category of overconvergent isocrystals. One embeds $X$ into a (good) admissible formal scheme $P$ which is smooth in the neighborhood of $X$. If we call $\bar X$ the Zariski closure of $X$ in $P$, an overconvergent isocrystal on $X \subset \bar X$ is a module with an integrable connection on a strict neighborhood of $X$ in $P$ whose Taylor series converges on a strict neighborhood of the diagonal. The point is that this definition is essentially independent of the choice of $P$ when $\bar X$ is fixed. Moreover, it is also independent of $\bar X$ if $\bar X$ is proper. This is the hard part of the theory. In our approach, the embedding of $X$ into $P$ defines a sieve $X_{P}/S$ of the topos of sheaves on $\textrm{AN}^\dagger(X/S)$. We will show that when $\bar X$ is proper, this is a covering sieve.
We are thus reduced to consider the localized category $\textrm{AN}^\dagger(X_{P}/S)$ with the induced sheaf $\mathcal O^\dagger_{X_{P}/S}$. We will show that the category of $\mathcal O^\dagger_{X_{P}/S}$-modules of finite presentation is equivalent to the category of overconvergent isocrystals on $X \subset \bar X/S$.
This is easier because we can use the bridge provided by modules with stratifications.

I want to finish this introduction with a description of our site $\textrm{AN}^\dagger(X/S)$.
We start with the following remark : any embedding $X \hookrightarrow P$ should define an object and our sheaves are defined on the the generic fibre $P_{K}$ of $P$. It is therefore natural to consider the pairs $(X \hookrightarrow P, P_{K} \hookleftarrow V)$ composed of a formal embedding and an open immersion.
Moreover, our sheaves should mainly depend on the tube of $X$. It is therefore necessary to identify tow objects with the same trace on the tube. Morphisms are just pairs of compatible morphisms and the topology is the analytic topology. This is the small overconvergent site of $X/S$. In this article, we will consider the big overconvergent site where we consider more general pairs $(U \hookrightarrow P, P_{K} \leftarrow V)$ where $U$ is a variety over $X$ and $V \to P_{K}$ any morphism.
This is better setting to study functoriality questions.

\section*{Acknowledgments}

I had the chance to talk to Vladimir Berkovich, Pierre Berthelot, Bruno Chiarellotto, Antoine Ducros, Michel Gros, Christian Naumovic, Michael Temkin and Nobuo Tsuzuki about questions related to this article and I want to thank them here.
Part of this work was done at the universities of Hiroshima, Nagoya and Tokyo in march 2004 where I was invited by Nobuo Tsuzuki, Kazuhiro Fujiwara and Atsushi Shiho, respectively.
The support of the European network Arithmetic Algebraic Geometry was also very valuable.

\section*{Warning}
In this article, valuations are non-trivial and affinoid algebras, affinoid spaces and analytic varieties are always strict. It is very likely that most results extend straightforwardly to the general situation but I have been to lazy to check the references.

\section*{Notations}
A \emph{presheaf} on a category $C$ is simply a contravariant functor to the category of sets.
We will denote by $\hat C$ the category of presheaves on $C$ and implicitly consider $C$ as a subcategory (representable presheaves).
Recall that a \emph{sieve} of an object $X$ of $C$ is simply a subobject of the corresponding presheaf.
A topology on $C$ is a family of sieves of objects of $C$, called \emph{covering sieves}, satisfying some properties.
When $C$ is endowed with a topology, we call it a \emph{site}.
A presheaf $F$ is a \emph{sheaf} if $F(X) = Hom(R, F)$ whenever $R$ is a covering sieve.
We will denote by $\tilde C$ the category of sheaves of sets on $C$.
This is a \emph{topos}.

\section*{Conventions}
Throughout this article, $K$ is a (non trivial) complete ultrametric field with valuation ring $\mathcal V$ and residue field $k$.

%
%
\section{Formal embeddings}

We denote by
$$
\mathcal V\{T_1, \ldots, T_n\} = \{\sum_{\underline i \geq 0} a_{\underline i} T^{\underline i}, a_{\underline i} \in \mathcal V, |a_{\underline i}| \to 0\}.
$$
the ring of convergent power series over $\mathcal V$.
A formal $\mathcal V$-scheme will always be assumed to have a locally finite covering by formal affine schemes $\textrm{Spf}(A)$ where $A$ is a quotient $\mathcal V\{T_1, \ldots, T_n\}$ by an ideal of finite type. They form a category $\textrm{FSch}(\mathcal V)$.
A formal scheme is said \emph{admissible} if it is $\mathcal V$-flat or, in other words, if it has no torsion.
If $S$ is a formal $\mathcal V$-scheme, we will also denote by $\textrm{FSch}(S)$ the category of formal schemes over $S$, which is just the localized category
$$
\textrm{FSch}(S) := \textrm{FSch}(\mathcal V)_{/S}.
$$ 
We will make a regular use of this notion of \emph{localized category}: if $C$ is a category and $T \in C$, the localized category $C_{/T}$ is the category of all morphisms $s : X \to T$ with $X \in C$.
Actually, this can be generalized to the case where $T$ is only a presheaf on $C$ and we consider sections $s \in T(X)$ and morphisms compatible with theses sections.
Of course, these two notions coincide when $T$ is representable.

\bigskip

If $K \hookrightarrow K'$ is an isometric embedding of complete ultrametric fields, and $\mathcal V'$ denotes the valuation ring of $K'$, there exists an extension functor
$$
\textrm{FSch}(\mathcal V) \to \textrm{FSch}(\mathcal V'), \quad P \mapsto P_{\mathcal V'}
$$
with $P_{\mathcal V'} = \textrm{Spf}(\mathcal V' \hat \otimes_{\mathcal V}ÊA)$ when $P = \textrm{Spf} (A)$.
It might be necessary to consider the category of \emph{generalized formal schemes} over $\mathcal V$.
An object is a pair $(P', K')$ where $K'$ is an isometric extension of $K$ and $P'$ a formal scheme on its valuation ring $\mathcal V'$.
Usual formal schemes correspond to the case $K' = K$.
A morphism $(P'', K'') \to (P', K')$ is made of an isometric $K$-embedding $K' \hookrightarrow K''$ and a morphism of $\mathcal V''$-formal schemes morphism $P'_{\mathcal V''} \to P''$.

If $X$ is any scheme, we denote by $\textrm{Sch}(X)$ the category of schemes locally of finite presentation over $X$.
If $k$ denotes the residue field of $K$, the category $\textrm{Sch}(k)$ may (and will) be seen as a full subcategory of $\textrm{FSch}(\mathcal V)$ and its object will be called algebraic varieties.
Moreover, the embedding $\textrm{Sch}(k) \hookrightarrow \textrm{FSch}(\mathcal V)$ has a right adjoint 
$$
\textrm{FSch}(\mathcal V) \to \textrm{Sch}(k), \quad P \mapsto P_k,
$$
sending a formal scheme to its special fibre.

A \emph{formal embedding} $X \hookrightarrow P$ (or $X \subset P$ for short) is a (locally closed) immersion over $\mathcal V$ of a $k$-scheme into a formal $\mathcal V$-scheme. A \emph{morphism of formal embeddings}
$$
(f \subset v) : (X' \subset P') \to (X \subset P)
$$
is a pair of morphisms (over $k$ and $\mathcal V$ respectively)
$$
(f : X' \to X, v : P' \to P),
$$
such that the diagram
$$
\begin{array} {ccc} X' & \hookrightarrow & P'\\ \;\;\;\downarrow f && \;\;\; \downarrow v\\
X & \hookrightarrow & P \end{array}
$$
is commutative.
When $f$ is the identity, we will just say that $v$ is a \emph{morphism of formal embeddings of $X$}.

\bigskip

\textbf{Example : } In order to connect our construction to Monsky-Washnitzer's, we may consider the following situation: We let $A$ be a $\mathcal V$-algebra of finite type;
The choice of a presentation of $A$ gives an embedding
$$
\textrm{Spec} A \hookrightarrow \mathbf A^N_{\mathcal V} \subset \mathbf P^N_{\mathcal V}
$$
which can be used to embed $X := \textrm{Spec} A_k$ into the formal completion $P := \widehat{\mathbf P^N_{\mathcal V}}$ of $\mathbf P^N_{\mathcal V}$.
This is a formal embedding.

\begin{prop} We have the following results :

\begin{enumerate}

\item With obvious composition, formal embeddings $(X \subset P)$ form a category $\textrm{Fmb}(\mathcal V)$ with finite inverse limits.

\item The forgetful functor
$$
\textrm{Fmb}(\mathcal V) \to \textrm{FSch}(\mathcal V), \quad (X \subset P) \mapsto P
$$
is exact and has an adjoint on the right
$$
\textrm{FSch}(\mathcal V) \to \textrm{Fmb}(\mathcal V), \quad Q \mapsto (Q_k \subset Q).
$$
\item The forgetful functor
$$
\textrm{Fmb}(\mathcal V) \to \textrm{Sch}(k), \quad (X \subset P) \mapsto X
$$
(is left exact and) has an adjoint on the left
$$
\textrm{FSch}(k) \to \textrm{Fmb}(\mathcal V), \quad Y \mapsto (Y \subset Y).
$$
\end{enumerate}

\end{prop}

\textbf {Proof : }
It should be clear that $\textrm{Fmb}(\mathcal V)$ is indeed a category. The first adjointness assertion simply says that, if $X \subset P$ is a formal embedding, any morphism of formal $\mathcal V$-schemes $P \to Q$ sends $X$ into $Q_k$. The other adjointness assertion is trivial. Finally, we have to check that, if we are given a finite diagram $(X_i \subset P_i)$, it has an inverse limit which is just $(\varprojlim X_i \subset \varprojlim P_i)$. But this again should be clear. $\Box$\bigskip

If $T$ is a presheaf on $\textrm{Fmb}(\mathcal V)$, we will call
$$
\textrm{Fmb}(T) := {\textrm{Fmb}(\mathcal V)}_{/T}
$$
the \emph{category of formal embeddings over $T$}.
This applies in particular to the case of a formal embedding $T := (X \subset P)$ over $\mathcal V$ (identified with the presheaf that it represents) in which case we obtain the category of formal embeddings over $X \subset P$. As a special case, if we still write $S$ for the image $(S_k \subset S)$ in $\textrm{Fmb}(\mathcal V)$ of a formal $\mathcal V$-scheme $S$, we get the category $\textrm{Fmb}(S)$ of formal embeddings over $S$.
An object of $\textrm{Fmb}(S)$ is simply a formal embedding $X \hookrightarrow P$ into a formal $S$-scheme and a morphism in $\textrm{Fmb}(S)$ is a morphism of formal embeddings $(f, u)$ where $u$ is an $S$-morphism.

\bigskip

We will finish with the description of a less trivial but fundamental case but we should recall first that any functor $g : C' \to C$ induces a functor on presheaves
$$
\hat g^{-1} : \hat C \to \hat C', \quad T \mapsto g \circ T
$$
which has a right adjoint $\hat g_{*}$ (and also a left adjoint $\hat g_{!}$).

\bigskip

If $S$ is a formal $\mathcal V$-scheme, then the forgetful functor $\textrm{Fmb}(S) \to Sch(S_k)$ induces a morphism $I_{S*}$ on presheaves.
On the other hand, the forgetful functor $\textrm{Fmb}(S) \to \textrm{Fmb}(\mathcal V)$ induces a morphism $J_{S*}$ on presheaves which has a left adjoint $J_{S}^{-1}$ which in turn has a left adjoint $J_{S!}$.
This is a standard result on localization.
If $X$ is an $S_{k}$-scheme, we may consider the presheaf $X/S := J_{S!}I_{S*}X$, which comes naturally with a morphism to $S$.
This presheaf is not representable in general but we can easily give a down to earth description of the category $\textrm{Fmb}(X/S)$ of formal embeddings over $X/S$ :

\begin{prop} Let $S$ be a formal $\mathcal V$-scheme and $X$ an $S_k$ scheme. A formal embedding over $X/S$ is an immersion $U \hookrightarrow P$ over $S$ of an $X$-scheme $U$ into a formal $S$-scheme.
A morphism of formal embeddings
$$
(f \subset v) : (U' \subset P') \to (U \subset P)
$$
is a pair of morphisms (over $X$ and $S$ respectively)
$$
(f : U' \to U, v : P' \to P),
$$
such that the diagram
$$
\begin{array} {ccc} U' & \hookrightarrow & P'\\ \;\;\;\downarrow f && \;\;\; \downarrow v\\
U & \hookrightarrow & P \end{array}
$$
is commutative.

\end{prop}

\textbf {Proof : } An object of $\textrm{Fmb}(X/S)$ is an object $(U \subset P)$ of $\textrm{Fmb}(\mathcal V)$ together with a section of $ j_{S!}I_{S*}X$ over $(U \subset P)$. By definition, such a section corresponds to a morphism $(U \subset P) \to (S_k \subset S)$ together with a section of $I_{S*}X$ on this object. Finally, by adjunction, such a section corresponds to a morphism $U \to X$ over $S_k$. The assertion on morphisms is proved in the same way. $\Box$\bigskip

As already mentioned, we will consider analytic spaces in the sense of Berkovich. We will denote by $\textrm{AN}(K)$ the category of analytic spaces over $K$, and more generally, if $V$ is an analytic space over $K$, we will denote by $\textrm{AN}(V) := \textrm{AN}(K)_{/V}$ the category of analytic spaces over $V$.
We will also have to consider the category of \emph{generalized analytic spaces} over $K$. An object is a pair $(V', K')$ where $K'$ is an isometric extension of $K$ and $V'$ an analytic space over $K'$.
Usual analytic spaces correspond to the case $K' = K$. A morphism $(V'', K'') \to (V', K')$ is made of an isometric $K$-embedding $K' \hookrightarrow K''$ and a morphism $V'_{K''} \to V''$.

As shown in \cite {Berkovich94}, section 1, there is a \emph{generic fiber} functor
$$
\textrm{FSch}(\mathcal V) \to \textrm{AN}(K), \quad P \mapsto P_K
$$
which is easily seen to be left exact.
When $P = \textrm{Spf} A$, we simply have $P_K = \mathcal M(A_K)$.
If $P$ is a formal $\mathcal V$-scheme, there is a natural \emph{specialization map}
$$
sp : P_K \to P, x \mapsto \tilde x.
$$
When $P = \textrm{Spf} A$, any $x \in P_K$ induces a continuous morphism $A_K \to \mathcal K(x)$ which reduces to a morphism $A_k \to k(x)$ whose kernel is $\tilde x \in S_k$.
Note that specialization is anticontinuous when $\textrm{FSch}(\mathcal V)$ is endowed with its Zariski topology and $\textrm{AN}(K)$ is endowed with its analytic topology.
More precisely, the inverse image of an open subset is a closed analytic domain and the inverse image of a closed subscheme is an open subset, as the following local description shows.
In fact, this will not be a problem for us for we will implicitly endow $\textrm{FSch}(\mathcal V)$ with the coarse topology (all presheaves are sheaves).
If $X \subset P$ is a formal embedding, we will consider the \emph{tube} $]X[_P := sp^{-1}(X)$ of $X$ in $P$. This definition will be generalized later. When $P = \textrm{Spf}(A)$ and
$$
X := \{x \in P, \forall i = 1, \ldots, r, f_i(x) = 0 \textrm{ and } \exists j = 1, \ldots s, g_j(x) \neq 0\},
$$
we have
$$
]X[_P := \{x \in P_K, \forall i = 1, \ldots, r, |f_i(x)| < 1 \textrm{ and } \exists j = 1, \ldots s, |g_j(x) = 1\}.
$$
Unlike in rigid geometry, it is not always true that an analytic space is locally affinoid. We will call a formal embedding $X \subset P$ \emph{good} if any point of $]X[_P$ has an affinoid neighborhood in $P_K$. We will also say that $P$ is \emph{good} at $X$.
Of course, saying that $P_K$ is good in Berkovich sense means that the trivial embedding of $P_k$ in $P$ is good. And in this case, any formal embedding into $P$ is good.
This is the case, for example if $P$ is affine or $P$ proper over $S$ and $S_K$ itself is good.
Recall however that if $P$ is the formal affine plane minus one point, then $P_K$ is not good (\cite {Temkin00}).
Nevertheless, if $X$ is the affine plane minus one point, then the embedding of $X$ into the formal affine plane $P$ is a good embedding.

%
%
\section{Geometric results} \label {tech}

The content of this section is not needed before the end of (the last) section \ref{comp}.
Therefore, the reader who wishes so can jump to the next section and come back later to this one. Note also that only very little material from section 1 is needed here

\bigskip

I thank V. Berkovich and M. Temkin for their help in understanding the geometric results of this section. Of course, they are not responsible for any mistake that may appear.

\bigskip

If $P$ is a formal $\mathcal V$-scheme and $X$ is a subset of $P_k$, we will denote by $\bar X^P$ or simply $\bar X$, the Zariski closure of $X$ in $P_k$.
A morphism of admissible formal schemes $u : P' \to P$ is said \emph{proper at $x' \in P'_k$} \label{properat} if the induced map $\overline {\{x'\}} \to P_k$ is proper.
It is said \emph{\'etale}, (resp. \emph{smooth}, resp. \emph{flat}, resp. \emph{proper}) at $X' \subset P'$ if it is \'etale (resp. smooth, resp. flat, resp. proper) at all $x' \in X'$.
Note that if $X'$ is a subscheme of $P'$, then $u$ is proper at $X'$ if and only if the restriction of $u$ to any irreducible component of $\bar X'$ is proper over $P_k$.
When $X'$ is quasi-compact, this just means that $\bar X'$ is proper over $P_k$.

\begin{prop} \label {Tem}
A morphism of admissible formal $\mathcal V$-schemes $u : P' \to P$ is proper at $X' \subset P$ if and only if $]X'[_{P'}$ is contained in the interior of $P'_K$ relative to $P_K$.
\end{prop} 

\textbf {Proof : }This easily follows from Theorem 4.1 of \cite {Temkin00} as explained in remark 5.8 of \cite {Temkin04}). $\Box$\bigskip

\begin{cor}
 Let $X \subset P$ be a good admissible formal embedding and
$$
(f, u) : (X' \subset P') \to (X \subset P)
$$
a morphism of admissible embeddings. If $u$ is proper at $X'$, then $P'$ is good at $X'$.
\end{cor}

\textbf {Proof : }It follows from Proposition \ref{Tem} that $u$ induces a morphism without boundary (i.e. closed in Berkovich sense) from a neighborhood $V'$ of $X'$ in $P$ to a neighborhood $V$ of $X$ in $P$. We may assume that $V$ is good and it follows that $V'$ is also good. $\Box$\bigskip

\begin{prop} \label{flat}
Let
$$
(f, u) : (X' \subset P') \to (X \subset P)
$$
be a morphism of good formal embeddings with $u$ flat at $X'$. Then, $u$ induces a universally flat morphism $V' \to V$ between good neighborhoods of $X'$ and $X$ in $P'_K$ and $P_K$, respectively.
\end{prop}

\textbf {Proof : } Clearly, $u$ induces a morphism $V' \to V$ between good neighborhoods of $X'$ and $X$ in $P'_K$ and $P_K$, respectively. Now, if $x' \in ]X'[_{P'}$, there exists affine neighborhoods $Q'$ of $\tilde x'$ in $P'$ and $Q$ of $u(\tilde x')$ in $P'$ such that the induced morphism $Q' \to Q$ is flat. Since flatness of formal schemes is stable under base change, it follows that the corresponding morphism $Q'_K \to Q_K$ is universally flat. Restriction to open subsets still gives a universally flat morphism
$$
V' \cap Q'_K \to V \cap Q_K
$$
of good analytic spaces. In particular, it is universally flat at $x'$. Now, since $V \cap Q_K \hookrightarrow V$ and $V' \cap Q'_K \hookrightarrow V'$ are inclusions of good analytic domains into good analytic varieties, they are universally flat. It follows that $u_K : V' \to V$ is universally flat at $x'$. Since this is true for all points of $]X'[_{P'}$ and that universal flatness is a local notion, we can shrink $V$ and $V'$ in order to get a universally flat morphism $V' \to V$. $\Box$\bigskip

\begin{lem} \label {pointwise}
Let $X \subset P$ be a good formal embedding. If $P$ is smooth (resp. \'etale) at $X$, there exists a quasi-smooth (resp. quasi-\'etale) good neighborhood $V$ of $]X[_P$ in $P_K$.
\end{lem}

\textbf {Proof : } We know that there exists a good neighborhood $V$ of $]X[_P$ in $P_K$. Let $x$ be a point of $]X[_P$ and $Q$ a smooth affine neighborhood of $\tilde x$ in $P$. Since $\Omega^1_Q$ is locally free (rep. $0$), so is $\Omega^1_{Q_K}$ and it follows from Proposition 6.23 of \cite {Ducros03*} that $Q_K$ is quasi-smooth (resp. quasi-\'etale). Since $V \cap Q_K$ is also good, Corollary 6.21 of \cite {Ducros03*} tells us that $V$ is quasi-smooth (resp. quasi-\'etale) at $x$. This is true for any point $x$ in $]X[_P$ and it follows that $V$ is quasi-smooth (resp. quasi-\'etale) at each $x \in ]X[_P$. Since quasi-smoothness (resp. quasi-\'etaleness) is a local notion, we may shrink $V$ a little bit in order to get a quasi-smooth (resp. quasi-\'etale) good neighborhood of $]X[_P$ in $P_K$. $\Box$\bigskip

\begin{prop} \label {quasismooth}
 Let
$$
(f, u) : (X' \subset P') \to (X \subset P)
$$
be a morphism of formal embeddings with $X \subset P$ good and $u$ smooth (resp. \'etale) at $X'$.
Then, $u$ induces a quasi-smooth (resp. quasi-\'etale) morphism $V' \to V$ between good neighborhoods of $X'$ and $X$ in $P'_K$ and $P_K$, respectively.
\end{prop}

\textbf {Proof : } It follows from \ref{flat} that $u$ induces a universally flat morphism $V' \to V$ between good neighborhoods of $X'$ and $X$ in $P'_K$ and $P_K$, respectively. Using Proposition 6.27 of \cite{Ducros03*}, we may assume that $P$ is reduced to a point, in which case, this is just the assertion of Lemma \ref{pointwise}. $\Box$\bigskip

\begin{cor} \label {smooth}
Let
$$
(f, u) : (X' \subset P') \to (X \subset P)
$$
be a morphism of formal embeddings with $X \subset P$ good admissible and $u$ proper and smooth (resp. \'etale) at $X'$. Then, $u$ induces a smooth (resp. \'etale) morphism $V' \to V$ between good neighborhoods of $X'$ and $X$ in $P'_K$ and $P_K$, respectively.
\end{cor}

\textbf {Proof : } Since a smooth (resp. \'etale) morphism is simply a quasi-smooth (resp. quasi-\'etale) morphism without boundary, this assertion follows from Propositions \ref{Tem} and \ref{smooth}. $\Box$\bigskip

\begin{prop} \label {poly}
Let $u : P' \to P$ be a morphism of formal embeddings of $X$.
If $u$ is smooth (resp. \'etale) at $X$, the fibers of the induced morphism $]X[_{P'} \to ]X[_P$ are strict open polydiscs (resp. $]X[_{P'} \simeq ]X[_P$).
\end{prop}

\textbf {Proof : }Let $x \in ]X[_P$ and $\mathcal V(x)$ be the ring of integers of $\mathcal K(x)$, the completed residue field at $x$. We may extend scalars by $\textrm{Spf} \mathcal V(x) \to P$ and therefore assume that $P = \textrm{Spf} \mathcal V$, $P_K = \mathcal M(K)$ and $X = \textrm{Spec} k$. We have to show that $]\tilde x[_{P'}$ is a closed polydisc. We may replace $P'$ by any open neighborhood of $\tilde x$ and therefore assume that there is an \'etale morphism $P' \to \hat {\mathbf A}^d_\mathcal V$ sending $x$ to $0$. It follows from Lemma 4.4 of \cite {Berkovich99} that this morphism induces an isomorphism 
$$
]\tilde x[_{P'} \simeq ]0[_{\hat {\mathbf A}^d_\mathcal V} = \mathbf B^d(0, 1^{-}).
$$ $\Box$\bigskip

\begin{lem} \label {cart}
Let $u : P' \to P$ be a morphism of formal embeddings of $X$. Assume that the induced map $\bar X^{P'} \to \bar X^P$ is separated. Then, there exists a neighborhood $V'$ of $]X[_{P'}$ in $P'_K$ such that for all $x \in ]X[_P$, we have
$$
u_K^{-1}(]\tilde x[_P) \cap V' = ]\tilde x[_{P'}.
$$
In particular,
$$
u_K^{-1}(]X[_P) \cap V' = ]X[_{P'}.
$$
\end{lem}

\textbf {Proof : }We may take $V' := ]\bar X^{P'}[_{P'}$. It is therefore sufficient to show that
$$
u^{-1}(\tilde x) \cap \bar X^{P'} = \{\tilde x\}.
$$
We are led to check that the dense open immersion
$$
X \hookrightarrow u^{-1}(X) \cap \bar X^{P'}
$$
is bijective. But the projection
$$
u^{-1}(X) \cap \bar X^{P'} \to X
$$
is separated as a pull-back of a separated map. A dense open immersion that admits a separated section is necessarily an isomorphism. $\Box$\bigskip

We will need the following lemma which is inspired by proposition 3.7.5 of \cite {Berkovich93}.

\begin{lem} \label {germs}
Let
$$
v : (V', x') \to (V, x)
$$
be a smooth morphism of germs of good analytic varieties. Any isomorphism of germs
$$
\varphi_x : (v^{-1}(x), x') \simeq (\mathbf A^d_{\mathcal \mathcal K(x)}, 0)
$$
over $\mathcal K(x)$ extends to an isomorphism of germs
$$
(V', x') \simeq (A^d_V, (0, x))
$$
over $(V, x)$.
\end{lem}

\textbf {Proof : } We may assume that $v$ comes from a morphism of affinoid algebras $A \to A'$ and that our isomorphism comes from a morphism
$$
\varphi_x^* : \mathcal K(x)\{T_1/R, \ldots, T_d/R\} \to \mathcal K(x) \hat \otimes_A A'.
$$
for $R$ big enough. Since the image of $\mathcal O_{V,x} \otimes_A A'$ is dense in $\mathcal K(x) \hat \otimes_A A'$, there exists $$
\varphi_1, \ldots, \varphi_d \in \mathcal O_{V,x} \otimes_A A'
$$
such that $\|\varphi_x^*(T_i) - \bar \varphi_i \| < R$.

After an automorphism of $\mathcal K(x)\{T_1/R, \ldots, T_d/R\}$, we may assume that $\varphi_x^*(T_i) = \bar \varphi_i$. Moreover, shrinking $V$ (and consequently $V'$) if necessary, we may also assume that $\varphi_1, \ldots, \varphi_d \in A'$ and we can consider the induced morphism $\varphi : V' \to \mathbf A^d_V$. By construction, it induces the inclusion $\varphi_x$ between the fibers at $x$. It follows from Lemma 3.7.7 of \cite {Berkovich93} that it is an isomorphism in a neighborhood of $x'$. $\Box$\bigskip

\begin{thm} \label {main}
Let $u : P' \to P$ be a morphism of good admissible formal embeddings of $X$ which is proper and smooth at $X$. Then, for all $x \in ]X[_P$, there exists a neighborhood $V$ of $x$ in $P_K$ and a section $s : V \to P'_K$ of $u_K$ on $V$ such that $\widetilde{s(y)} = \tilde y$ for all $y \in V \cap ]X[_P$.
\end{thm}

\textbf {Proof : } 
We know from corollary \ref {smooth} that $u$ induces a smooth morphism $v : V' \to V$ between good neighborhoods of $X$ in $P'_K$ and $P_K$.
Moreover, we may assume thanks to lemma \ref {cart} that $v^{-1}(]\tilde y[_P) = ]\tilde y[_{P'}$ for $y \in ]X[_P$. In particular, any local section of $v$ will satisfy the last condition. 

Since $v^{-1}(]X[_P) = ]X[_{P'}$, it follows from proposition \ref{poly} that there exists an isomorphism
$$
\varphi_x : v^{-1}(x) \simeq \mathbf B^d_{\mathcal K(x)}(0, 1^{-}).
$$
We set $x' := \varphi_x^{-1}(0)$.
It follows from lemma \ref {germs} that there exists a neighborhood $W'$ of $x'$ in $V'$ such that $\phi_x$ extends to some open immersion $\varphi : W' \hookrightarrow \mathbf A^d_V$. Since $v$ is a flat morphism without boundary, the image of $W'$ in $V$ is an open neighborhood $W$ of $0$ (unpublished result of Berkovich, see lemma 6.2.8 of \cite{Ducros03*}) and we can take the zero section of $\mathbf A^d_W$. $\Box$\bigskip

%
%
\section{Analytic varieties}

An \emph{analytic variety over $\mathcal V$} is a couple made of a good formal embedding $X \subset P$ and a morphism of analytic varieties $\lambda : V \to P_K$.
It can be represented by the diagram
$$
X \hookrightarrow P \stackrel {sp}\leftarrow P_K \stackrel \lambda \leftarrow V
$$
The \emph{tube} of $X$ in $V$ is $]X[_V := \lambda^{-1}(]X[_P)$ and we will denote by
$$
i_{X,V} : ]X[_V \hookrightarrow V
$$
the inclusion map.
We will generally write $(X \subset P \leftarrow V)$ or $(X, V)$ to make notations shorter. We might also forget $\lambda$ or $sp$ in the notations and just write $sp : V \to P$ or $\lambda : V \to P$, in which case we will write $]X[_V = sp^{-1}(X)$ or $]X[_V = \lambda^{-1}(X)$.
Also, if $x \in ]X[_V$, we will write $\tilde x = sp(\lambda (x))$.
Finally, when no confusion should arise, we will simply call $i_X$, $i_V$ or even $i$ the inclusion map

\bigskip

\textbf{Example : } We can consider again the Monsky-Washnitzer situation.
We saw that if $A$ is a finitely presented $\mathcal V$-algebra, we can build a formal embedding
$$
X := \textrm{Spec} A_k \subset P := \widehat{\mathbf P^N_{\mathcal V}}.
$$
We may also consider the inclusion morphism
$$
V := (\textrm{Spec} A_K)^{an} \hookrightarrow (\mathbf{P}^N_K)^{an} = (\widehat{\mathbf{P}^N_{\mathcal V}})_K = P_K.
$$
in order to get an analytic variety
$$
X \subset P \stackrel {sp}\leftarrow P_K \hookleftarrow V
$$
over $\mathcal V$.
We have
$$
]X[_P = \mathbf B^N(0,1^+) \cap V.
$$
Actually, for each $\lambda > 1$, we may set
$$
V_\lambda := \mathbf B^N(0, \lambda) \cap V
$$
and we get another analytic variety
$$
X \subset P \stackrel {sp}\leftarrow P_K \hookleftarrow V_\lambda
$$
over $\mathcal V$.
Note that $V_\lambda$ is affinoid so that we can write $V_\lambda = \mathcal M(A_\lambda)$. Note also that $]X[_P = \mathcal M(\hat A_K)$, that $A_\lambda \subset \hat A_K$ for $\lambda$ close to $1$ and that $\cup A_\lambda \subset \hat A_K$ is the generic fiber $A^\dagger_K$ of the weak completion $A^\dagger$ of $A$.

\bigskip

We now come to the definition of morphisms. It should be remarked before that if $K \hookrightarrow K'$ is an isometric embedding and if $\mathcal V'$ and $k'$ denote the valuation ring and residue field of $K'$ as usual, then any analytic variety $(X \subset P \leftarrow V)$ over $\mathcal V$ gives rise to an analytic variety $(X_{k'} \subset P_{\mathcal V'} \leftarrow V_{K'})$ over $\mathcal V'$. 
A \emph{hard morphism of analytic varieties} from $(X', V')$ to $(X,V)$ is a couple of morphisms
$$
(f : X' \to X, u : V' \to V)
$$
such that
$$
\forall x \in ]X'[_{V'}, \widetilde {u(x)} = f(\widetilde {x})
$$
and the same condition holds after any isometric extension of $K$.
Note that it is sufficient to check the condition only for rational points (after any isometric extension of $K$).

In practice, we will write
$$
(f, u) : (X', V') \to (X, V).
$$
Most of the time, our hard morphism will fit into a commutative diagram
$$
\begin{array} {ccccccc}
X' & \hookrightarrow & P' & \stackrel {sp}\leftarrow & P'_K & \stackrel {\lambda'} \leftarrow & V' \\
\downarrow f & & \downarrow v & & \downarrow v_K && \downarrow u\\
X & \hookrightarrow & P & \stackrel {sp}\leftarrow & P_K & \stackrel \lambda \leftarrow & V
\end{array}.
$$
We will then say that it \emph{extends to formal schemes}.

\bigskip

\textbf{Example :} Back again to our Monsky-Washnitzer situation. We assume now that we are given two algebras of finite type $A$ and $A'$ over $\mathcal V$ and we use the same notations as above. A morphism $f : X' \to X$ corresponds to an algebra homomorphism $A_k \to A'_k$. If $A$ is smooth, this homomorphism lifts to a homomorphism $A^\dagger \to A'^\dagger$ which induces $A_\lambda \to A'_\mu$ for $\lambda$ and $\mu$ close enough to $1$ and therefore gives a morphism $u : V'_\mu \to V_\lambda$. It is clear that we get a hard morphism
$$
(f, u) : (X', V'_\mu) \to (X, V_\lambda)
$$
which in general will \emph{not} extend to formal schemes.

\begin{prop}
We have the following results :
\begin{enumerate}

\item With obvious composition of morphisms, analytic varieties over $\mathcal V$ and hard morphisms form a category ${\textrm{AN}(\mathcal V)}$ with finite inverse limits.

\item The forgetful functor
$$
{\textrm{AN}}(\mathcal V) \to \textrm{AN}(K), \quad (X \subset P \leftarrow V) \mapsto V
$$
is exact with right adjoint
$$
\textrm{AN}(K) \to{\textrm{AN}}(\mathcal V), \quad W \mapsto (\textrm{Spec} k \subset \textrm{Spf} \mathcal V \leftarrow W).
$$
\item The forgetful functor
$$
{\textrm{AN}}(\mathcal V) \to \textrm{Sch}(k), \quad (X \subset P \leftarrow V) \mapsto X
$$
(is left exact and) has a left adjoint
$$
\textrm{Sch}(k) \to{\textrm{AN}}(\mathcal V), \quad Y \mapsto (Y = Y \leftarrow \emptyset).
$$
\item The functor
$$
\textrm{Fmb}(\mathcal V) \to{\textrm{AN}}(\mathcal V), \quad (X \subset P) \mapsto (X \subset P \leftarrow P_K)
$$
is left exact.

\end{enumerate}

\end{prop}

\textbf {Proof : } It is not difficult to see that the inverse limit of a diagram
$$
X_i \hookrightarrow P_i \stackrel {sp}\leftarrow P_{iK} \stackrel {\lambda_i}\leftarrow V_i
$$
indexed by some finite set $I$ is simply
$$
\varprojlim X_i \hookrightarrow \prod P_i \stackrel {sp}\leftarrow \prod P_{iK} \stackrel \lambda \leftarrow \varprojlim V_i.
$$
Moreover, if all the morphisms extend in a compatible way to formal schemes, the limit is
$$
\varprojlim X_i \hookrightarrow \varprojlim P_i \stackrel {sp}\leftarrow \varprojlim P_{iK} \stackrel \lambda \leftarrow \varprojlim V_i.
$$
All our assertions easily follow from these remarks. $\Box$\bigskip

It is also important to remark that the assignment
$$
(X \subset P \leftarrow V) \leadsto (X \subset P)
$$
is not functorial because of our very loose definition of morphisms in ${\textrm{AN}}(\mathcal V)$: they do not always extend to formal schemes.

\bigskip

A hard morphism
$$
(f, u) : (U', V') \to (U, V).
$$
will induce a morphism
$$
]f[_u : ]X'[_{V'} \to ]X[_V
$$
between the tubes giving rise to a functor
$$
{\textrm{AN}}(\mathcal V) \to \textrm{AN}(K), \quad (X, V) \mapsto ]X[_{V}.
$$
Since $]f[_u$ is just the morphism induced by $u$, we will sometimes write $u : ]X'[_{V'} \to ]X[_V$.
Also, when $u = \textrm{Id}_V$, we will write $]f[_V : ]X'[_{V} \to ]X[_V$.

\begin{prop}
We have the following results :

\begin{enumerate}

\item Any hard morphism of analytic varieties over $\mathcal V$ is the composition of a morphism of the form
$$
(f, \textrm{Id}_{V}) : (X' \subset P' \leftarrow V) \to (X \subset P \leftarrow V)
$$
and a morphism of the form
$$
(\textrm{Id}_X, u) : (X \subset P \leftarrow V') \to (X \subset P \stackrel {\lambda} \leftarrow V).
$$
\item A hard morphism $(f, u) : (X', V') \to (X, V)$ is a hard isomorphism if and only if $f$ is an isomorphism, $u$ is an isomorphism and the induced map $]f[_{u} : ]X'[_{V'} \to ]X[_V$ is surjective (and therefore also an isomorphism).

\item Up to hard isomorphism, any hard morphism of analytic varieties over $\mathcal V$ extends to formal schemes.
More precisely, given any $(X', V') \to (X, V)$, there exists a hard isomorphism $(X'', V'') \simeq (X', V')$ such that the composite $(X'', V'') \to (X, V)$ extends to formal schemes.

\end{enumerate}

\end{prop}

\textbf {Proof : } The first assertion is clear : the morphism
$$
(f, u) : (X' \subset P' \stackrel {\lambda'} \leftarrow V') \to (X \subset P \stackrel {\lambda} \leftarrow V)
$$
splits as the composition of
$$
(f, \textrm{Id}_{V'}) : (X' \subset P' \stackrel {\lambda'} \leftarrow V') \to (X \subset P \stackrel {\lambda \circ u} \leftarrow V')
$$
followed by
$$
(\textrm{Id}_X, u) : (X \subset P \stackrel {\lambda \circ u} \leftarrow V') \to (X \subset P \stackrel {\lambda} \leftarrow V).
$$
Now, if $(f, u)$ is an isomorphism, then $f$ and $u$ must be isomorphisms themselves. Conversely, if $f$ and $u$ are both isomorphism, they induce an embedding of analytic domains
$$
]f[_u : ]X'[_{V'} \hookrightarrow ]X[_V
$$
It is clear that $f^{-1}$ and $u^{-1}$ define a hard morphism if and only if this embedding is an isomorphism.

Finally, any hard morphism
$$
(f, u) : (X' \subset P' \leftarrow V') \to (X \subset P \leftarrow V)
$$
gives rise to a commutative diagram
$$
\begin{array} {ccccccc}
X' & \hookrightarrow & P' & \stackrel {sp}\leftarrow & P'_K & \stackrel {\lambda'} \leftarrow & V' \\
|| & & \uparrow p_2& & \uparrow p_{2K} && ||\\
X' & \hookrightarrow & P \times P' & \stackrel {sp}\leftarrow & P_K \times P'_K & \stackrel {(\lambda \circ u, \lambda')} \leftarrow & V' \\
\downarrow f & & \downarrow p_1 & & \downarrow p_{1K} && \downarrow u \\
X & \hookrightarrow & P & \stackrel {sp}\leftarrow & P_K & \stackrel \lambda \leftarrow & V
\end{array}
$$
and it follows from the previous assertion that the upper morphism is a hard isomorphism. $\Box$\bigskip

\begin{cor} If
$$
X \hookrightarrow P \stackrel {sp}\leftarrow P_K \stackrel \lambda \leftarrow V
$$
is any analytic variety over $\mathcal V$ such that $\lambda$ factors through $P'_K$ where $P'$ is a formal subscheme of $P$ containing $X$, we get a hard isomorphism
$$
(X \subset P' \leftarrow V) \simeq (X \subset P \leftarrow V).
$$
\end{cor} $\Box$\bigskip

\textbf{Example : }
In the Monsky-Washnitzer situation, we can replace $P := \widehat {\mathbf P^N_\mathcal V}$ by the completion $P'$ of the algebraic closure of $\textrm{Spec} A$ in $\mathbf P^N_\mathcal V$ and get an isomorphic analytic variety over $X/\mathcal V$. We see that $X$ is open in $P'_k$ and $V$ is open in $P'_K$. This is a more pleasant situation to work with.

\begin{cor} If $(X \subset P \leftarrow V)$ is an analytic variety over $\mathcal V$ and $P' \to P$ is a formal blowing up centered outside $X$, then $P'_{K} = P_{K}$ and the induced hard morphism
$$
(X \subset P' \leftarrow V) \to (X \subset P \leftarrow V)
$$
is an isomorphism.
\end{cor} $\Box$\bigskip

\begin{cor} Any analytic variety over $\mathcal V$ is isomorphic to some analytic variety $(X \subset P \leftarrow V)$ where, if we denote as usual by $\bar X$ the Zariski closure of $X$ in $P$, $\bar X \backslash X$ is a divisor in $\bar X$.
\end{cor} $\Box$\bigskip

\begin{prop} \label {limitpro}
The functor
$$
{\textrm{AN}}(\mathcal V) \to \textrm{AN}(K), \quad (X \subset P \leftarrow V) \mapsto ]X[_{V}
$$
is left exact
\end{prop}

\textbf {Proof : } Since pull back is left exact, it is sufficient to prove the following :
if we are given a finite diagram $\{X_i\}_{i \in I}$ and formal embeddings $X_i \subset P_i$, then 
$$
]\varprojlim X_i[_{\prod P_i} = \varprojlim] X_i[_{P_i}.
$$
We only need to consider two cases. We assume first that we are given a finite set of embeddings $\{X_i \subset P_i\}_{i \in I}$, and we check that
$$
]\prod X_i[_{\prod P_i} = \prod ]X_i[_{P_i}.
$$
The second step is to show that if we are given a finite family of formal embeddings $\{X_i \subset P\}$, then
$$
]\cap X_i[_{P} = \cap] X_i[_{P}.
$$
Both results are standard and can be proved using the local description of the tubes. $\Box$\bigskip

A functor $g : C' \to C$ between tow sites is said \emph{cocontinuous} if $\hat g_{*}$ preserves sheaves.

\bigskip

The \emph{analytic topology} on ${\textrm{AN}}(\mathcal V)$ is the coarsest topology making cocontinuous the forgetful functor
$$
{\textrm{AN}}(\mathcal V) \to \textrm{AN}(K), \quad (X, V) \mapsto V.
$$
\begin{prop}
The analytic topology on ${\textrm{AN}}(\mathcal V)$ is generated by the following pretopology : families
$$
\{(X \subset P \leftarrow V_i) \to (X \subset P \leftarrow V)\}_{i \in I}
$$
that extend to the identity on $P$ and where $V = \cup_{i} V_i$ is an open covering.
\end{prop}

\textbf {Proof : } Note first that such families do define a pretopology. Now, the forgetful functor is cocontinuous if whenever $\{V_i \to V\}_{i \in I}$ is a covering, the family of all
$$
(X' \subset P' \leftarrow V') \to (X \subset P \leftarrow V)
$$
such that $V' \to V$ factors through some $V_i$, is a covering.
This condition is equivalent to
$$
(X' \subset P' \leftarrow V') \to (X \subset P \leftarrow V)
$$
factoring through some $(X \subset P \leftarrow V_i)$.
Thus, we see that the forgetful functor is cocontinuous if whenever $\{V_i \to V\}_{i \in I}$ is a covering, so is
$$
\{(X \subset P \leftarrow V_i) \to (X \subset P \leftarrow V)\}_{i \in I}. \quad \Box
$$
\bigskip

Recall that the \emph{canonical topology} on a category is the finest topology for which representable presheaves are sheaves.
A \emph{standard site} is a site with finite inverse limits where the topology is coarser than the canonical topology.

\begin{cor}
The site ${\textrm{AN}}(\mathcal V)$ is a standard site.
\end{cor}

\textbf {Proof : }
We have to show that the analytic topology is coarser than the canonical topology on ${\textrm{AN}}(\mathcal V)$.
We take any $(X \subset P \leftarrow V) \in {\textrm{AN}}(\mathcal V)$ and we have to prove that the presheaf
$$
(X' \subset P' \leftarrow V') \mapsto \textrm{Hom}((X' \subset P' \leftarrow V'), (X \subset P \leftarrow V))
$$
is a sheaf.
Thus, we are given a covering
$$
V' = \cup_{i} V'_{i}
$$
and a compatible family of morphisms
$$
\{(f, u_{i}) : (X' \subset P' \leftarrow V'_{i}) \to (X \subset P \leftarrow V)\}_{i \in I}.
$$
As Berkovich showed, the analytic topology is coarser than the canonical topology on $\textrm{AN}(K)$, and it follows that the $u_i's$ glue to a morphism $u : V' \to V$.
It is clear that $(f, u)$ is a morphism. $\Box$\bigskip

Note that if $P$ is the disjoint union of two copies of $\textrm{Spf} \mathcal V$ and if we embed $\textrm{Spec} k$ into $P$ on the left and send $\mathcal M(K)$ to $P_K$ on the right, the unique morphism
$$
(\textrm{Spec} k \subset P \leftarrow \mathcal M(K) \to (\textrm{Spec} k \subset \textrm{Spf} V \leftarrow \mathcal M(K))
$$
is not a covering although it is the identity both on the left and on the right.

\bigskip

Recall that a functor $g : C' \to C$ between two sites is said \emph{continuous} if $\hat g^{-1}$ preserves sheaves.
If the extension of $g$ to sheaves is exact, it defines a morphism of sites in the other direction $C \to C'$.
When the topology of $C$ is coarser than the canonical topology and $C'$ has finite inverse limits, it is actually sufficient to check that $g$ itself is left exact.
When $g$ is cocontinuous, the situation is nicer :
it automatically gives rise to a morphism of toposes $\tilde C \to \tilde C'$.

\bigskip

Unless otherwise specified, categories of schemes and formal schemes are always endowed with the coarse topology - and not the Zariski topology.

\begin{prop}

\begin{enumerate}

\item The functor
$$
\textrm{Fmb}(\mathcal V) \to{\textrm{AN}}(\mathcal V), \quad (X \subset P) \to (X \subset P \leftarrow P_K).
$$
is left exact and continuous, giving rise to a morphism of sites
$$
{\textrm{AN}}(\mathcal V)Ê\to \textrm{Fmb}(\mathcal V).
$$
\item The forgetful functor
$$
{\textrm{AN}}(\mathcal V) \to \textrm{AN}(K), \quad (X, V) \mapsto V
$$
is left exact, continuous and cocontinuous, giving rise to a morphism of sites
$$
\textrm{AN}(K) \to{\textrm{AN}}(\mathcal V).
$$
and a morphism of toposes
$$
\widetilde{{\textrm{AN}}(\mathcal V)} \to \widetilde{\textrm{AN}(K)}.
$$
\item The forgetful functor
$$
{\textrm{AN}}(\mathcal V) \to \textrm{Sch}(k), \quad (X, V) \mapsto X
$$
is left exact and continuous, giving rise to a morphism of sites $\textrm{Sch}(k) \to{\textrm{AN}}(\mathcal V)$. Moreover, the corresponding functor
$$
I_{\mathcal V*} : \widetilde{\textrm{Sch}(k)} \to \widetilde{\textrm{AN}(\mathcal V)}
$$
is fully faithful.

\item The tube functor
$$
{\textrm{AN}}(\mathcal V) \to \textrm{AN}(K), \quad (X, V) \mapsto ]X[_V
$$
is left exact and continuous, giving rise to a morphism of sites $\textrm{AN}(K) \to \textrm{AN}(\mathcal V)$.

\end{enumerate}

\end{prop}

\textbf {Proof : } The first assertion follows from the definition. And the second one is an immediate consequence of our hypothesis. 

Concerning the third assertion, exactness is already known and continuity is easily checked.
It is therefore sufficient to notice the full faithfulness of the stupid functor
$$
Y \mapsto (Y = Y \leftarrow \emptyset).
$$
For the last assertion, there is nothing to do. $\Box$\bigskip

%
%
\section{The overconvergent site}

Let $(X \subset P \leftarrow V)$ be an analytic variety over $\mathcal V$ and $W$ an analytic domain in $V$.
If $W$ is a neighborhood of $]X[_V$, the morphism
$$
(X \subset P \leftarrow W) \hookrightarrow (X \subset P \leftarrow V)
$$
is called a \emph{strict neighborhood}.
It is said \emph{open} if $W$ is open in $V$.
Sometimes, we might just say that $W$ is a strict neighborhood of $X$ in $V$.
It simply means that $W$ is a neighborhood of $]X[_{V}$ in $V$.
The relation with Berthelot's definition of strict neighborhood is highlighted by the following result.

\begin{prop} \label{strictneigh}
Let $(X \subset P \leftarrow V)$ be an analytic variety over $\mathcal V$ and $W$ an analytic domain in $V$.
Denote by $\bar X$ the Zariski closure of $X$ in $P$.
Then, $W$ is a strict neighborhood of $X$ in $V$ if and only if the covering 
$$
]\bar X[_V = ]\bar X[_W \cup ]\bar X \backslash X[_V
$$
is admissible.
\end{prop}

\textbf{Proof : } Assume first that $W$ is a strict neighborhood of $X$ in $V$. Replacing $W$ by some open neighborhood, we may assume that $W$ is open in $V$. Moreover, since $]\bar X[_V$ is open in $V$, we may also assume that $V = ]\bar X[_V$. In this case, the above covering is simply an open covering and we are done.
Conversely, if the covering is admissible and $x \in ]X[_V$, then $x \in ]\bar X[_W$ but $x \not \in ]\bar X \backslash X[_V$. It follows that $]\bar X[_W$ is a neighborhood of $x$ in $]\bar X[_V$, and \textsl{a fortiori} $W$ is a neighborhood of $x$ in $V$. $\Box$\bigskip

\begin{lem}
We have the following results :

\begin{enumerate}

\item A strict neighborhood is a monomorphism in ${\textrm{AN}}(\mathcal V)$.

\item Any composition of strict neighborhoods is a strict neighborhood.

\item Any pull back of a strict neighborhood is a strict neighborhood.

\end{enumerate}

\end{lem}

\textbf {Proof : } Since the forgetful functor
$$
{\textrm{AN}}(\mathcal V) \to \textrm{Sch}(k) \times \textrm{AN}(K)
$$
is obviously faithful, the first assertion follows from the fact that the identity $\textrm{Id}_X$ as well as the inclusion $W \hookrightarrow V$ are both monomorphisms. The second assertion is trivial. Finally, one easily checks that the pull-back of
$$
(X, W) \hookrightarrow (X,V)
$$
by
$$
(f,u) : (X', V') \to (X, V)
$$
is simply
$$
(X', u^{-1}(W)) \hookrightarrow (X',V'). \quad \Box
$$
\bigskip

\begin{prop}
The category ${\textrm{AN}}(\mathcal V)$ admits calculus of right fractions with respect to strict neighborhoods.
\end{prop}

\textbf {Proof : }Since strict neighborhoods form a subcategory made of monomorphisms which is stable by pull-backs, all the conditions of (\cite{GabrielZisman67}, I, 2.2.2) are satisfied. $\Box$\bigskip

The quotient category $\textrm{AN}^\dagger(\mathcal V)$ is the \emph{category of analytic varieties over $\mathcal V$}.
The description of $\textrm{AN}^\dagger(\mathcal V)$ is quite simple.
First, it has the same objects as ${\textrm{AN}}(\mathcal V)$.
Second, if $(X, V)$ and $(X', V')$ are two analytic varieties, then
$$
\textrm{Hom}_{\textrm{AN}^\dagger(\mathcal V)}((X', V'), (X, V)) = \varinjlim {\textrm{Hom}}_{{\textrm{AN}}(\mathcal V)}((X', W'), (X, V))
$$
where $W'$ runs through all neighborhoods of $]X[_V$ in $V$.
In other words, we may always replace $V$ by any neighborhood of $]X[_V$ in $V$ and get an isomorphic object.
From this down to earth description, it is not difficult to see that, like ${\textrm{AN}(\mathcal V)}$, the category $\textrm{AN}^\dagger(\mathcal V)$ has finite inverse limits.

\bigskip

\textbf{Example : } Back again to the Monsky-Washnitzer situation, we see that any homomorphism $A^\dagger \to A'^\dagger$ gives rise to a morphism $(X', V') \to (X, V)$.
It is not clear to me wether the converse is true or not.

\bigskip

If $C'$ is a site and $g : C' \to C$ any functor, the \emph{image topology} on $C$ si the coarsest topology that makes $g$ continuous.

\begin{lem} \label{catlem}
Let $C$ be a site with fibered products, $Q$ a set of morphisms in $C$ and $C_{Q}$ the quotient category.

\begin{enumerate} \renewcommand{\labelenumi}{\roman{enumi})}

\item If $C_Q$ is endowed with the image topology, the canonical map $C \to C_Q$ extends to an inverse image for an embedding of sites $u : C_Q \hookrightarrow C$.

\item The topology of $C_Q$ is generated by covering sieves of the form $\hat u^{-1}R \subset X' \to X$ where $R \subset X'$ is a covering sieve and $X' \to X$ is in $Q$.

\item The functor $u_*$ induces an equivalence between $\tilde C_Q$ and the full subcategory of $T \in \tilde C$ such that $T(\varphi)$ is an isomorphism whenever $\varphi \in Q$.

\end{enumerate}
\end{lem}

\textbf {Proof : } It is clear that $u^{-1}$ is left exact and, by definition, it is continuous. Moreover, by definition again, the composition functor $\hat u_* : \hat C_Q \to \hat C$ identifies $\hat C_Q$ with the full subcategory of $T \in \hat C$ such that $T(\varphi)$ is an isomorphism if $\varphi \in Q$. The first assertion follows.
Since $\hat u^{-1}$ is left exact, it is a general fact that $u^{-1}$ is continuous if and only if whenever $R$ is a sieve of $X \in C$, then $\hat u^{-1}R$ is a sieve of $X \in C_Q$. This proves the second assertion.
Now, we have to show that under our new assumptions, a presheaf $T$ on $C_Q$ is a sheaf if and only if $\hat u_* T$ is a sheaf. The condition is always necessary. By adjonction, it is also sufficient. Namely, we always have
$$
\textrm{Hom}(\hat u^{-1}R, T) = \textrm{Hom}(R, \hat u_* T) = \textrm{Hom}(X, \hat u_* T) = \textrm{Hom}(X, T). \quad \Box
$$
\bigskip

The image of the analytic topology on ${\textrm{AN}}(\mathcal V)$ is the \emph{analytic topology} on $\textrm{AN}^\dagger(\mathcal V)$.
The associated topos will be writen $\mathcal V_{AN^\dagger}$.

\begin{prop}
The canonical functor
$$
{\textrm{AN}}(\mathcal V) \to \textrm{AN}^\dagger(\mathcal V)
$$
is the inverse image for an embedding of sites
$$
\textrm{AN}^\dagger(\mathcal V) \hookrightarrow{\textrm{AN}}(\mathcal V).
$$
This embedding induces an equivalence between $\mathcal V_{AN^\dagger}$ and the full subcategory of sheaves $\mathcal F$ on ${\textrm{AN}}(\mathcal V)$ such that $\mathcal F(X, V) = \mathcal F(X, W)$ whenever $W$ is a strict neighborhood of $X$ in $V$.
\end{prop}

\textbf {Proof : } Follows from lemma \ref{catlem}. $\Box$\bigskip

Note that the analytic topology is \emph{compatible with strict neighborhoods} in the following sense : If $(X, V)$ is an analytic variety, $V = \cup_{i \in I} V_{i}$ is an open covering and for each $i \in I$, $W_i$ is a strict neighborhood of $X$ in $V_i$, there exists strict neighborhoods $W'_{i}$ of $X$ in $W_{i}$ such that $W' = \cup_{i \in I} W'_i$ is an open covering and $W'$ is a strict neighborhood of $X$ in $V$.

\bigskip
\textbf{Remark} : The Grothendieck topology is not compatible with strict neighborhoods as the following example shows :
$$
V = \mathbf P^1, \quad X = \emptyset,
$$
$$
W_1 = \mathbf D(0, 1^-) \subset V_1 = \mathbf D(0, 1^+), \quad W_2 = V_2 := \mathbf D(\infty, 1^+).
$$
\begin{prop} The analytic topology on $\textrm{AN}^\dagger(\mathcal V)$ is defined by the following pretopology : families of hard morphisms
$$
\{(X \subset P \leftarrow V_{i}) \to (X \subset P \leftarrow V)\}_{i \in I}
$$
where $\{V_{i}\}_{i \in I}$ is a an open covering of a strict neighborhood of $X$ in $V$.
\end{prop}

\textbf {Proof : } It follows from part 2 of lemma \ref{catlem} above that our topology is generated by these families and it is therefore sufficient to show that this is a pretopology.
Our families clearly contain the identities and stability by pull back is easily verified. It only remains to prove transitivity. This follows from the fact that the analytic topology is compatible with strict neighborhoods. $\Box$\bigskip

\begin{prop} The site $\textrm{AN}^\dagger(\mathcal V)$ is a standard site.
\end{prop}

\textbf {Proof : } As already mentioned, our site has fibered products and it remains to show that the topology is coarser than the canonical topology.
We let $(X \subset P \leftarrow V) \in \textrm{AN}(\mathcal V)$ and we show that the presheaf
$$
(X' \subset P' \leftarrow V') \mapsto \textrm{Hom}((X' \subset P' \leftarrow V'), (X \subset P \leftarrow V))
$$
is a sheaf on $\textrm{AN}^\dagger(\mathcal V)$.
So assume that we are given an analytic covering
$$
\{(X' \subset P' \leftarrow V'_{i}) \to (X' \subset P' \leftarrow V')\}_{i \in I}
$$
and a compatible family of morphisms
$$
\{(X' \subset P' \leftarrow V'_{i}) \to (X \subset P \leftarrow V)\}_{i \in I}.
$$
There exists strict neighborhoods $W'_i$ of $X'$ in $V'_i$ such that this last family is actually given by hard morphisms
$$
\{(X' \subset P' \leftarrow W'_{i}) \to (X \subset P \leftarrow V)\}_{i \in I}.
$$
Since the analytic topology is compatible to strict neighborhoods, we may assume that we have an open covering $W' = \cup_{i} W'_i$ and $W'$ is a strict neighborhood of $X'$ in $V$.
Now, everything is defined in ${\textrm{AN}}(\mathcal V)$ and we are done. $\Box$\bigskip

We recall again that categories of schemes and formal schemes are endowed with the coarse topology (which makes specialization continuous).

\begin{prop}

\begin{enumerate}

\item The functor
$$
\textrm{Fmb}(\mathcal V) \to \textrm{AN}^\dagger(\mathcal V), \quad (X \subset P) \to (X \subset P \leftarrow P_K)
$$
is left exact and continuous, giving rise to a morphism of sites
$$
\textrm{AN}^\dagger(\mathcal V)Ê\to \textrm{Fmb}(\mathcal V).
$$
\item The forgetful functor
$$
\textrm{AN}^\dagger(\mathcal V) \to \textrm{Sch}(k), \quad (X, V) \mapsto X
$$
is left exact and continuous, giving rise to a morphism of sites
$$
I^\dagger_{\mathcal V} : \textrm{Sch}(k) \to \textrm{AN}^\dagger(\mathcal V).
$$
\item The tube functor
$$
\textrm{AN}^\dagger(\mathcal V) \to \textrm{AN}(K), \quad (X, V) \mapsto ]X[_V
$$
is left exact and continuous, giving rise to a morphism of sites $\textrm{AN}(K) \to \textrm{AN}^\dagger(\mathcal V)$.

\end{enumerate}

\end{prop}

\textbf {Proof : } The first assertion directly follows from the analogous results for ${\textrm{AN}}(\mathcal V)$.
And the second one is almost trivial.
Finally, the last assertion also follows from the analogous result for ${\textrm{AN}}(\mathcal V)$. $\Box$\bigskip

We should also remark that the assignment
$$
(X \subset P \leftarrow V) \leadsto V
$$
is no longer functorial.

\begin{prop} 
Let $K \hookrightarrow K'$ be an isometric embedding and $\mathcal V'$, $k'$ denote the valuation ring and residue field of $K'$ respectively. Then the extension functor 
$$
\textrm{AN}^\dagger(\mathcal V) \to \textrm{AN}^\dagger(\mathcal V'), \quad (X \subset P \leftarrow V) \mapsto (X_{k'} \subset P_{\mathcal V'} \leftarrow V_{K'})
$$
is left exact and continuous, giving rise to a morphism of sites
$$
\textrm{AN}^\dagger(\mathcal V') \to \textrm{AN}^\dagger(\mathcal V) 
$$
\end{prop}

\textbf {Proof : } Both assertions directly follow from our definition. $\Box$\bigskip

Note that the analogous statement for $\textrm{AN}$ is also true but the analytic site was only a tool for us and we should not have to mention it in the future.

Recall that if $C$ is a site and $g : C' \to C$ any functor, then the \emph{induced topology} on $C'$ is the finest topology that makes $g$ continuous.
Note that when $g$ is left exact and $C'$ has fibered products, a family in $C'$ is a covering family for the induced topology if and only if its image in $C$ is a covering family.
This applies in particular to the case of a localization functor $C_{/T} \to C$ where $C$ is a site and $T$ a presheaf on $C$.

\bigskip

An \emph{overconvergent presheaf} on $\mathcal V$ is a presheaf $T$ on $\textrm{AN}^\dagger(\mathcal V)$.
We will call $\textrm{AN}^\dagger(T) := \textrm{AN}^\dagger(\mathcal V)_{/T}$ the \emph{category of analytic varieties over $T$}.
The category $\textrm{AN}^\dagger(T)$ is endowed, as explained above, with the analytic topology induced by the analytic topology on $\textrm{AN}^\dagger(\mathcal V)$.
We get the \emph{overconvergent site $\textrm{AN}^\dagger(T)$ of $T$} and the corresponding overconvergent topos $T_{\textrm{AN}^\dagger}$ whose objects are \emph{overconvergent sheaves}Êon $T$.
Thus, we see that in general, an analytic variety over $T$ is a triple $(U, V, t)$ where $(U, V)$ is an analytic variety over $\mathcal V$ and $t \in T(U, V)$. A morphism $(U', V', t') \to (U, V, t)$ being a morphism of analytic varieties $(f, u) : (U', V') \to (U, V)$ such that $T(f, u)(t') = t$.
Of course any morphism of presheaves $f : T' \to T$ induces a morphism of toposes 
$$
f : \widetilde{\textrm{AN}^\dagger(T')} \to \widetilde{\textrm{AN}^\dagger(T)}.
$$
As a first example of overconvergent presheaf, we can consider the case of an analytic variety $(X, V)$ identified with the corresponding presheaf.
We obtain the category of analytic varieties over $(X, V)$.
An object is simply a morphism $(X', V') \to (X, V)$ and a morphism in $\textrm{AN}^\dagger(X, V)$ is simply a usual morphism that commutes with the given ones.
Of course, any morphism of analytic varieties $(f, u) : (X', V') \to (X, V)$ will induce a morphism of toposes
$$
f : \widetilde{\textrm{AN}^\dagger(X', V')} \to \widetilde{\textrm{AN}^\dagger(X, V)}.
$$
As a particular case, if $S$ is a formal $\mathcal V$-scheme, we will call
$$
\textrm{AN}^\dagger(S) := \textrm{AN}^\dagger(S_k, S_K)
$$
the \emph{category of analytic varieties over $S$}. 
By functoriality, any morphism of formal $\mathcal V$-schemes $v : S' \to S$ provides us with a morphism of toposes
$$
f : \widetilde{\textrm{AN}^\dagger(S')} \to \widetilde{\textrm{AN}^\dagger(S)}.
$$
\begin{prop} Let $S$ be a formal $\mathcal V$-scheme.
Then, up to isomorphism, an analytic variety over $S$ is is a couple made of a good formal embedding $X \subset P$ over $S$ and a morphism of analytic varieties $\lambda : V \to P_K$ over $S_K$.
And a morphism is just a morphism $(f, u) : (X', V') \to (X, V)$ of analytic varieties over $\mathcal V$ where $f$ is an $S_k$-morphism and $u$ an $S_K$-morphism.
\end{prop}

\textbf{Proof : } We know that, up to isomorphism, any morphism of analytic varieties extends to formal schemes. More precisely, if $(X \hookrightarrow P \leftarrow V)$ is an analytic variety over $S$, we can embed $X$ in $P_S$ and $V$ in $(P_S)_K$ in order to get another analytic variety $(X \hookrightarrow P_S \leftarrow V)$ over $S$ which is clearly isomorphic to the original one. Everything else follows. $\Box$\bigskip

If $S$ is a formal $\mathcal V$-scheme, we have $I_{\mathcal V}^{^\dagger-1}(S) = S_k$ and, the morphism of sites
$$
I^\dagger_\mathcal V : \textrm{AN}(k) \to \textrm{AN}^\dagger(\mathcal V)
$$
induces by localization, a morphism of sites
$$
I^\dagger_S : \textrm{AN}(S_k) \to \textrm{AN}^\dagger(S).
$$
Recall also that if
$$
J_{S} : \textrm{AN}^\dagger(S) \to AN^\dagger(\mathcal V)
$$
denotes the localization map, then $J_{S!}$ factors through an isomorphism
$$
 \widetilde{\textrm{AN}^\dagger(S)} \simeq \widetilde{\textrm{AN}^\dagger(\mathcal V)}_{/\tilde S}.
$$
Any $S_k$ scheme $X$ defines a sheaf on $S_k$ and we let $X/S := J_{S!} I^\dagger_{S*}X$.
Thus, if $S$ is a formal $\mathcal V$-scheme and $X$ is an $S_k$-scheme, we obtain the important category $\textrm{AN}^\dagger(X/S)$ of analytic varieties over $X/S$.
Note that $X/S$ comes naturally with a morphism to $S$ in $\widetilde{AN^\dagger(S)}$.

\begin{prop}
Let $S$ be a formal $\mathcal V$-scheme and $X$ an $S_k$-scheme. Then, up to isomorphism, an analytic variety over $X/S$ is is a couple made of a good formal embedding $U \subset P$ over $S$ with $U$ an $X$-scheme and a morphism of analytic varieties $\lambda : V \to P_K$ over $S_K$.
And a morphism is just a morphism $(f, u) : (X', V') \to (X, V)$ of analytic varieties over $S$ where $f$ is an $X$-morphism.
\end{prop}

\textbf{Proof : } By definition, an analytic variety over $X/S$ is given by an object $(U \subset P \leftarrow V)$ of $\textrm{AN}^\dagger(\mathcal V)$ and a section of the sheaf $J_{S!}I_{S*}X$ on this object.
Such a section si given by some structural morphism
$$
(U \subset P \leftarrow V) \to (S_k \subset S \leftarrow S_K)
$$
and a section of $I_{S*}X$ on $(U \subset P \leftarrow V)$ over $S$.
By adjunction, this section corresponds to a morphism from $U := \hat I_{S}^{-1}(U \subset P \leftarrow V)$ to $X$ over $S_k$. Summarizing, we get a good formal embedding $U \subset P$ over $S$ with $U$ an $X$-scheme and a morphism of analytic varieties $\lambda : V \to P_K$ over $S_K$ as expected. The last assertion is easily checked. $\Box$\bigskip

Any $S_k$-morphism $f : X' \to X$ induces, by functoriality, a morphism of toposes
$$
f : \widetilde{\textrm{AN}^\dagger(X'/S)} \to \widetilde{\textrm{AN}^\dagger(X/S)}.
$$
Also, given any analytic variety $(X, V)$ over $S$, we have a canonical morphism of toposes
$$
\widetilde{\textrm{AN}^\dagger(X,V)} \to \widetilde{\textrm{AN}^\dagger(X/S)}.
$$
\begin{prop} Let $\sigma : K \hookrightarrow K'$ be an isometry of complete ultrametric fields and denote as usual by $\mathcal V'$ and $k'$ the valuation ring and residue fields of $K'$.
Let $S'$ be a formal $\mathcal V'$-scheme and $X'$ be an $S'_{k'}$-scheme. Assume that we are given a morphism of formal schemes $v : S' \to S$ over $\mathcal V \to \mathcal V'$ and a morphism $f : X' \to X$ above $v$.
Then, there is a natural morphism of toposes
$$
\widetilde{\textrm{AN}^\dagger(X'/S')} \to \widetilde{\textrm{AN}^\dagger(X/S)}.
$$
\end{prop}

\textbf {Proof : }We first consider the composite morphism of toposes
$$
\widetilde{\textrm{AN}^\dagger(S')} \to \widetilde{\textrm{AN}^\dagger(S \hat\otimes \mathcal V')} \to \widetilde{\textrm{AN}^\dagger(S)}.
$$
The pull back of $X$ is simply $X \otimes_k k'$ and our morphism is induced by the canonical map $X' \to X \otimes_k k'$. $\Box$\bigskip

We can say a little more about our functor $I^\dagger$ :

\begin{prop}
The functor
$$
I^\dagger_{\mathcal V*} :\textrm{Sch}(k) \to \widetilde{\textrm{AN}^\dagger(\mathcal V)}
$$
is continuous and left exact, giving rise to a morphism of topos
$$
U^\dagger_{\mathcal V} : \widetilde{\textrm{AN}^\dagger(\mathcal V)} \to \widehat{\textrm{Sch}(k)}.
$$
And we have a sequence of adjoint functors
$$
{I^\dagger_{\mathcal V}}^{-1} \quad , \quad I^\dagger_{\mathcal V*} = {U^\dagger_{\mathcal V}}^{-1} \quad , \quad U^\dagger_{\mathcal V*}.
$$
with $U^\dagger_{\mathcal V} \circ I^\dagger_{\mathcal V} = \textrm{Id}$.
\end{prop}

\textbf {Proof : } The functor $I^\dagger_{\mathcal V*}$ being a direct image of morphism of sites is automatically left exact.
Moreover, it is also automatically continuous since we have the coarse topology on the left hand side. $\Box$\bigskip

Actually, localization gives us for each formal scheme $S$ and each $S_k$-scheme $X$, two morphisms of topos
$$
I^\dagger_{X/S} : \widetilde{\textrm{Sch}(X)} \to \widetilde{\textrm{AN}^\dagger(X/S)}.
$$
and
$$
U^\dagger_{X/S} : \widetilde{\textrm{AN}^\dagger(X/S)} \to \widetilde{\textrm{Sch}(X)}
$$
such that, at the topos level, we have a sequence of adjoint functors
$$
{I^\dagger_{X/S}}^{-1} \quad , \quad I^\dagger_{X/S*} = {U^\dagger_{X/S}}^{-1} \quad , \quad U^\dagger_{X/S*},
$$
and $U^\dagger_S \circ I^\dagger_S = \textrm{Id}$ again.
Finally, one can check that
$$
U^\dagger_{S*}(\mathcal F)(X) = \Gamma(X/S, \mathcal F).
$$

\bigskip

We finish this section with a description of inverse image with respect to an immersion of algebraic varieties.

\begin{prop} \label{rep}
Let $S$ be a formal $\mathcal V$-scheme and $\alpha: Y \hookrightarrow X$ be an immersion of algebraic varieties over $S_{lk}$.
Then, if $(U, V) \in \textrm{AN}^\dagger(X/S)$, we have $\hat \alpha^{-1}(U, V) = (p^{-1}(Y), V)$ where $p : U \to X$ is the structural map.
\end{prop}

\textbf {Proof : }By definition, if $(U', V') \in \textrm{AN}^\dagger(Y/S)$, a section of $\hat \alpha^{-1}(U, V)$ on $(U', V')$ is a morphism
$$
(f, u) : (U', V') \to (U, V)
$$
in $\textrm{AN}^\dagger(X/S)$.

Then $f$ factors over $X$ through some morphism $f' : U' \to p^{-1}(Y)$ and we clearly get a morphism
$$
(f', u) : (U', V') \to (p^{-1}(Y), V)
$$
in $\textrm{AN}^\dagger(Y/S)$ since $p^{-1}(Y) \subset U$.
Conversely, any morphism in $\textrm{AN}^\dagger(Y/S)$ gives by composition such a morphism in $\textrm{AN}^\dagger(X/S)$. $\Box$\bigskip

%
%
\section{Realization of sheaves}

Let $V$ be a topological space and $T$ a subset of $V$.
An inclusion $V'' \subset V'$ of open subsets of $V$ is called a \emph{$T$-isomorphism} if $V' \cap T = V'' \cap T$. A family $V'_i$ of open subsets of $V$ is called an open \emph{$T$-covering} of an open subset $V'$ if for each $i$, we have $V'_i \cap T \subset V'$ and also $V' \cap T \subset \cup_{i} V'_i$.

\begin{lem} \label {relative}
Let $V$ be a topological space and $T$ a subset of $V$.
Then
\begin{enumerate}
\item The category $\textrm{Open}(V)$ of open subsets of $V$ admits calculus of right fractions with respect to $T$-isomorphisms.
\item The image topology on the quotient category $\textrm{Open}(V)_{T}$ is generated by the pretopology of open $T$-coverings.
\item The inclusion $T \subset V$ induces an equivalence of sites
$$
\textrm{Open}(V)_{T} \simeq \textrm{Open}(T)
$$
when $T$ is equipped with the induced topology.
\end{enumerate}
\end{lem}

\textbf {Proof : }It is clear that $T$-isomorphisms form a subcategory which is stable by intersection with an open subset. The first assertion formally follows. Now, the inclusion $T \subset V$ clearly induces an equivalence of categories
$$
\textrm{Open}(V)_{T} \simeq \textrm{Open}(T).
$$
And open $T$-coverings in $\textrm{Open}(V)$ correspond to open coverings in $T$. $\Box$\bigskip

As usual, any analytic variety $V$ over $K$ will be endowed with its analytic topology and we will denote by $V_{\textrm{an}}$ the small topos of sheaves on $V$.

\begin{prop}
Let $(X, V)$ be an analytic variety over $\mathcal V$.
The obvious functor $V' \mapsto (X, V')$ on $\textrm{Open}(V)$ induces a functor
$$
\textrm{Open}(V)_{]X[_V} \to \textrm{AN}^\dagger(X,V).
$$
which is continuous, cocontinuous and left exact.
\end{prop}

\textbf {Proof : }It is clear that if $V' \subset V''$ is an $]X[_{V}$-isomorphism, then $(X,V' )\hookrightarrow (X, V)$ is an open strict neighborhood. Moreover, it follows directly from our definitions that the induced functor is continuous and cocontinuous : any open $]X[_V$-covering of some open subset $V' \subset V$ gives rise to an analytic open covering of $(X, V')$ and conversely. It is also clearly left exact. $\Box$\bigskip

\begin{cor}
If $(X, V)$ is an analytic variety over $\mathcal V$, the functor $V' \mapsto (X, V')$ defines a morphism of sites
$$
\varphi_{X,V} : \textrm{AN}^\dagger(X,V) \to \textrm{Open}(V)_{]X[_V} \simeq \textrm{Open}(]X[_V).
$$
and a morphism of toposes
$$
\psi_{X,V} : (]X[_V)_{\textrm{an}} \to (X,V)_{\textrm{AN}^\dagger}
$$
giving rise to a sequence of adjoint functors at the topos level :
$$
\varphi^{-1}_{X,V} \quad , \quad \varphi_{X,V*} = \psi_{X,V}^{-1} \quad , \quad \psi_{X,V*}
$$
\end{cor} $\Box$\bigskip

\begin{prop}
If $(f, u) : (X', V') \to (X, V)$ is a morphism of analytic varieties over $\mathcal V$ and $\mathcal F$ a sheaf on $]X[_V$, the diagram
$$
\begin {array} {ccc}
(X',V')_{\textrm{AN}^\dagger} & \stackrel {\varphi_{X',V'}} \longrightarrow & (]X'[_{V'})_{\textrm{an}}\\
\downarrow j && \downarrow ]f[_{u} \\
(X,V)_{\textrm{AN}^\dagger} & \stackrel {\varphi_{X,V}} \longrightarrow & (]X[_V)_{\textrm{an}}
\end{array}
$$
is commutative.
\end{prop}

\textbf {Proof : }This follows from the fact that, if $W$ is an open subset of $V$, then $(f, u)^{-1}(X, W) = (X', u^{-1}(W))$. $\Box$\bigskip

\begin{cor} If $(f, u) : (X', V') \to (X, V)$ is a morphism of analytic varieties over $\mathcal V$ and $\mathcal F$ an overconvergent sheaf on $(X,V)$, one has
$$
\Gamma((X', V'), \varphi_{X,V}^{-1}\mathcal F) = \Gamma(]X'[_{V'}, ]f[_{u}^{-1}\mathcal F).
$$
\end{cor} $\Box$\bigskip

\begin{prop}
If $(X, V)$ is an analytic variety over $\mathcal V$, we have
$$
\varphi_{X,V} \circ \psi_{X,V} = Id.
$$
\end{prop}

\textbf {Proof : }We have
$$
(\varphi_{X,V} \circ \psi_{X,V})^{-1}(\mathcal F) = \psi_{X,V}^{-1}(\varphi_{X,V}^{-1} \mathcal F)
= \varphi_{X,V*}(\varphi_{X,V}^{-1} \mathcal F).
$$
Thus, if $V'$ is an open subset of $V$, we see that
$$
\Gamma(]X[_{V'}, (\varphi_{X,V} \circ \psi_{X,V})^{-1}(\mathcal F)) = \Gamma(]X[_{V'}, \varphi_{X,V*}(\varphi_{X,V}^{-1} \mathcal F))
$$
$$
= \Gamma((X, V'), \varphi_{X,V}^{-1} \mathcal F)= \Gamma(]X'[_{V'}, \mathcal F).
$$ $\Box$\bigskip

If $S$ is a formal $\mathcal V$-scheme, we can apply the above considerations to $S$, identified with the analytic variety $(S_k, S_K)$ as usual. We obtain two morphism of toposes 
$$
\psi_{S} : S_{K\textrm{an}} \to S_{\textrm{AN}^\dagger}
$$
and
$$
\varphi_{S} : S_{\textrm{AN}^\dagger} \to S_{K\textrm{an}}.
$$
By composition, we get that for any overconvergent presheaf $T$ over $S$, there is a canonical morphism of toposes
$$
p_{T} : T_{\textrm{AN}^\dagger} \to S_{\textrm{AN}^\dagger} \to S_{K\textrm{an}}.
$$
We now turn to the definition of the realizations.
If $(X, V)$ is an analytic variety over $\mathcal V$ and $\mathcal F$ is an overconvergent sheaf on $(X,V)$, then the sheaf
$$
\mathcal F_{X,V} :=  \varphi_{X,V*} \mathcal F
$$
on $]X[_V$ will be called the \emph{realization} of $\mathcal F$ on $]X[_V$.
We will often write $\mathcal F_V$ instead of $\mathcal F_{X,V}$. In the case $V = P_K$, we will write $\mathcal F_{X \subset P}$ or $\mathcal F_P$ and call it the \emph{realization} of $\mathcal F$ on $X \subset P$ or $P$.
By extension, if $T$ is an overconvergent presheaf over $\mathcal V$, $(X, V)$ an analytic variety over $T$ and $\mathcal F$ an analytic sheaf on $T$, one defines the \emph{realization} $\mathcal F_{X,V}$ of $\mathcal F$ on $(X,V)$ as follows : We first take the inverse image of $\mathcal F$ by the localization morphism
$$
(X,V)_{\textrm{AN}^\dagger} \to T_{\textrm{AN}^\dagger}
$$
and then, its realization on $(X, V)$.
There is a simpler way to define the realizations of $\mathcal F$:
Note first that, since the topology of $]X[_V$ is induced by the topology of $V$, any open subset of $]X[_V$ is of the form $]X[_{V'}$ with $V'$ open in $V$. Moreover, $\mathcal F(X, V')$ is independent of $V'$. Then, we simply have
$$
\mathcal F_{X,V}(]X[_{V'}) = \mathcal F(X, V').
$$
Using realizations, it is not difficult to describe the morphism $\varphi_{X,V}$. Of course, by definition, for any overconvergent sheaf $\mathcal F$ on $(X, V)$, we have $\varphi_{X,V*} \mathcal F = \mathcal F_{X,V}$.
But also, if $\mathcal F$ is a sheaf on $]X[_{V}$ and $(f,u) : (X', V') \to (X, V)$, then
$$
(\varphi_{X,V}^{-1}\mathcal F)_{X',V'} = ]f[_{u}^{-1}\mathcal F.
$$
Finally, if $T$ is an overconvergent presheaf on $\mathcal V$, $(f, u) : (X', V') \to (X, V)$ a morphism of analytic varieties over $T$ and $\mathcal F \in T_{\textrm{AN}}$, functoriality gives us a morphism
$$
\phi_{fu} : ]f[_u^{-1} \mathcal F_{X,V} \to \mathcal F_{X, V'}
$$
on $]X'[_{V'}$. Of course the morphisms $\phi_{fu}$ satisfy the usual compatibility condition. These data uniquely determine $\mathcal F$ as the following proposition states.

\begin{prop}
If $T$ is an overconvergent presheaf on $\mathcal V$, the category $T_{\textrm{AN}}$ is equivalent to the following category :
\begin{enumerate}
\item An object is a collection of sheaves $\mathcal F_{X,V}$ on $]X[_V$ for each $(X, V) \in \textrm{AN}^\dagger(T)$ and morphisms $\phi_{fu\mathcal F} : ]f[_u^{-1} \mathcal F_{X,V} \to \mathcal F_{X, V'}$ for each $(f, u) : (X', V') \to (X, V)$, satisfying the usual compatibility condition.
\item A morphism is a collection of morphisms $\mathcal F_{X,V} \to \mathcal G_{X,V}$ compatible with the morphisms $\phi_{fu}$.
\end{enumerate}
\end{prop}

\textbf {Proof : } This is completely standard. $\Box$\bigskip

\begin{prop}
If $T$ is an overconvergent presheaf on $\mathcal V$, the topos $T_{\textrm{AN}^\dagger}$ has enough points. More precisely, if $(X, V) \in \textrm{AN}^\dagger(T)$ and $x \in ]X[_V$, then the functor $\mathcal F \mapsto \mathcal F_{V,x}$ is a fiber functor and they form a conservative family.
\end{prop}

\textbf {Proof : }The point $x \in ]X[_V$ defines a point of the topos $]X[_{V\textrm{an}}$ and composition with $\psi_{X,V}$ followed by the structural morphism $(X,V)_{\textrm{AN}} \to T_{\textrm{AN}^\dagger}$ gives a point of our topos. We have to show that this family of points is conservative. Since the points of a topological space form a conservative family, it is sufficient to note that the ``family'' of realization functors $\{\mathcal F \mapsto \mathcal F_{X,V}\}_{(X, V) \in \textrm{AN}^\dagger(T)}$ is faithful. $\Box$\bigskip

Again, we consider the case of an immersion of algebraic varieties.

\begin{prop} \label {immer}
Let $S$ be a formal $\mathcal V$-scheme, $\alpha : Y \hookrightarrow X$ an immersion of $S_k$-schemes, and $(U, V) \in \textrm{AN}^\dagger(X/S)$. Let $p : U \to X$ be the structural map, $U' = p^{-1}(Y)$ and $]\alpha[ : ]U'[_V \hookrightarrow ]U[_V$ the inclusion map. If $\mathcal F$ is a sheaf on $\textrm{AN}^\dagger(X/S)$, then
$$
(\alpha_{\textrm{AN}^\dagger*}\mathcal F)_V = ]\alpha[_{*}\mathcal F_V.
$$
\end{prop}

\textbf {Proof : }We have seen in proposition \ref{rep} that $\alpha^{-1}_{\textrm{AN}}(U, V) = (U', V)$, and it follows that $(\alpha_{\textrm{AN}*}\mathcal F)(U, V) = \mathcal F(U', V)$. And the same is true for any open subset of $V$. Thus we see that
$$
i_{U*}(\alpha_{\textrm{AN}^\dagger*}\mathcal F)_V = i_{U'*}\mathcal F_V = i_{U*}]\alpha[_{*}\mathcal F_V.
$$
Pulling back by $i_U$ gives the result. $\Box$\bigskip

\begin{cor}
If $\alpha : Y \hookrightarrow X$ is an open immersion, then $\alpha_{\textrm{AN}^{\dagger}*}$ is exact.
\end{cor}

\textbf {Proof : }In this case $]\alpha[_{*}$ is the inclusion of a closed subset. $\Box$\bigskip

%
%
\section{Coherent sheaves}

We prove here some general results about coherent sheaves on locally compact spaces and we apply them to Berkovich analytic spaces.

\begin{prop} \label {limit2}
Let $(V, \mathcal O_V)$ be a ringed space and $i : T \hookrightarrow V$ the inclusion of a closed subspace. Assume $V$ paracompact or $T$ compact.
If $\mathcal F$ is an $\mathcal O_V$-module of finite presentation and $\mathcal G$ any $\mathcal O_V$-module, the canonical map
$$
\varinjlim \textrm{Hom}_{\mathcal O_{V'}} (\mathcal F_{|V'}, \mathcal G_{|V'}) \to \textrm{Hom}_{i^{-1} \mathcal O_V}(i^{-1} \mathcal F, i^{-1} \mathcal G),
$$
where $V'$ runs through all open neighborhoods of $T$, is an isomorphism.
\end{prop} 

\textbf {Proof : }Since $\mathcal F$ is finitely presented, we have
$$
i^{-1}\mathcal H\textrm{om}_{\mathcal O_V} (\mathcal F, \mathcal G) = \mathcal H\textrm{om}_{ i^{-1}\mathcal O_V} ( i^{-1}\mathcal F, i^{-1}\mathcal G)
$$
and therefore, 
$$
\textrm{Hom}_{i^{-1} \mathcal O_V}(i^{-1} \mathcal F, i^{-1} \mathcal G) = \Gamma (T, i^{-1} \mathcal H\textrm{om}_{\mathcal O_V} (\mathcal F, \mathcal G)).
$$
Since $T$ is compact, or closed in $V$ paracompact, it follows from \cite {KashiwaraSchapira90}, Proposition 2.5, that
$$
 \Gamma (T, i^{-1} \mathcal \mathcal H\textrm{om}_{\mathcal O_V} (\mathcal F, \mathcal G)) = \varinjlim \Gamma (V', \mathcal H\textrm{om}_{\mathcal O_V} (\mathcal F, \mathcal G)_{|V'})
$$
and we know that local hom commute to localization. $\Box$\bigskip

\begin{cor} \label {limit3} With the same hypothesis, if $\mathcal F$ and $\mathcal G$ are two $\mathcal O_V$-modules of finite presentation such that $i^{-1} \mathcal F = i^{-1} \mathcal G$, there exists an open neighborhood $V'$ of $T$ in $V$ such that $\mathcal F_{|V'} = \mathcal G_{|V'}$.
\end{cor} $\Box$\bigskip 

\begin{prop} \label{icoher}
Let $(V, \mathcal O_V)$ be a locally compact ringed space and $i : T \hookrightarrow V$ the inclusion of a locally closed subspace.
If $\mathcal O_V$ is a coherent ring, then $i^{-1}\mathcal O_V$ is also a coherent ring.
\end{prop}

{\bf Proof :} We may clearly assume that $T$ is closed in $V$.
We have to show that if $W$ is an open subset of $V$ and $m : i^{-1}\mathcal O_{W}^N \to i^{-1}\mathcal O_{W}$ a morphism, then $\ker m$ is of finite type.
Given any $x \in T$, there exists an open subset $\mathcal U$ of $V$ and a compact subset $K$ of $V$ such that $x \in \mathcal U \subset K \subset W$.
Since $K \cap T$ is compact, it follows from Proposition \ref {limit2} that there exists an open neighborhood $\mathcal U'$ of $K \cap T$ in $W$ and a morphism $n : \mathcal O_{\mathcal U'}^N \to \mathcal O_{\mathcal U'}$ whose restriction to $K \cap T$ coincides with the restriction of $m$.
The same thing is obviously true over $\mathcal U \cap T$. Since $\mathcal O_V$ is coherent, the kernel of $n$ is of finite type. Since pull back is exact, the restriction of $m$ to $\mathcal U \cap T$ is also of finite type. And we are done. $\Box$\bigskip

\begin{prop}
Let $(V, \mathcal O_V)$ be a locally compact ringed space and $i : T \hookrightarrow V$ the inclusion of a closed subspace. Assume that $T$ is compact or $V$ is countable at infinity. If $\mathcal F$ is a coherent $i^{-1}\mathcal O_V$-module of finite presentation, there exists an open neighborhood $V'$ of $T$ in $V$ and a coherent $\mathcal O_{V'}$-module $\mathcal G$ such that $\mathcal F = i_{V'}^{-1}\mathcal G$.
\end{prop}

\textbf {Proof : }Let $\mathcal F$ be a coherent $i^{-1}\mathcal O_V$-module.
For each $x \in T$, there exists a sequence $x \in \mathcal U \subset K \subset W$ with $\mathcal U, W$ open in $V$ and $K$ compact and a morphism $m : i^{-1}\mathcal O_{W}^{r} \to i^{-1}\mathcal O_{W}^{s}$ with $\mathcal F_{|W \cap T} = \textrm{coker} m$.
Since $K \cap T$ is compact, it follows from Proposition \ref{limit2} that there exists an open neighborhood $\mathcal U'$ of $K \cap T$ in $V$ and a morphism $n : \mathcal O_{\mathcal U'}^{r} \to \mathcal O_{\mathcal U'}^{s}$ whose restriction to $K \cap T$ coincides with the restriction of $m$. Let $\mathcal G' := \textrm{coker} n$. We may replace $\mathcal U$ with $\mathcal U \cap \mathcal U'$ in order to have $\mathcal U \subset \mathcal U'$. If we set $\mathcal G := \mathcal G'_{|\mathcal U}$, we have $i^{-1}\mathcal G = \mathcal F_{|\mathcal U \cap T}$.

Thus, if we assume that $T$ is compact, we see that we can find a finite open covering $\{\mathcal U_k\}_{k=1}^n$ of $T$ in $V$, and for each $k$ a coherent $\mathcal O_{\mathcal U_k}$ module $\mathcal G_k$ such that $i^{-1}\mathcal G_k = \mathcal F_{|\mathcal U_k \cap T}$. More precisely, there exists for each $k$ an open subset $\mathcal U'_k$ of $V$ with $\mathcal U_k \subset \mathcal U'_k$, a compact subset $K_k$ of $V$ with $K_k Ê\cap T \subset \mathcal U'_k$ and a coherent module $\mathcal G'_k$ on $\mathcal U'_k$ whose restriction to $K_k \cap T$ coincides with the restriction of $\mathcal F$. In particular, the restrictions of $\mathcal G'_k$ and $\mathcal G'_l$ to $K_k \cap K_l \cap T$ coincide. It follows from Corollary \ref{limit3} that there exists an open neighborhood $\mathcal U'_{kl}$ of $K_k \cap K_l \cap T$ in $V$ such that $(\mathcal G'_k)_{|\mathcal U'_{kl}} = (\mathcal G'_l)_{|\mathcal U'_{kl}}$.
Now, since there are only a finite number of them, we may shrink each $\mathcal U'_k$ in order to have for each pair $(k, l)$, $\mathcal U'_k \cap \mathcal U'_l \subset \mathcal U'_{kl}$. We can then, for each $k$, replace $\mathcal U_k$ by $\mathcal U_k \cap \mathcal U'_k$ and still get a covering of $T$. We have for each pair $(k,l)$, $(\mathcal G_k)_{|\mathcal U_{k} \cap \mathcal U_{l}} = (\mathcal G_l)_{|\mathcal U_{k} \cap \mathcal U_{l}}$. It follows that the $\mathcal G_k$ glue together in order to give a coherent $\mathcal O_{V'}$-module $\mathcal G$ such that $\mathcal F = i_{V'}^{-1}\mathcal G$.

We now consider the second case, namely we assume that $V$ is countable at infinity. In other words, we assume that $V$ is an increasing union of compact subsets $K_n, n \in \mathbf N$. Since, for each $n \in \mathbf N$, $T \cap K_n$ is compact in $K_n$, it follows from the first case that the restriction of our coherent $i^{-1}\mathcal O_V$-module $\mathcal F$ to $K_n \cap T$ extends to some coherent module $\mathcal G_n$ on some open neighborhood $V_n$ of $K_n \subset T$ in $V$. By induction, we may assume that $V_n \subset V_{n+1}$ and, using Corollary \ref{limit3} again, that $\mathcal G_{n+1|V_n} = \mathcal G_n$ so that they glue in order to give a coherent module $\mathcal G$ on $V' = \cup V_n$. $\Box$\bigskip

If $\mathcal O_V$ is a sheaf of rings on some topological space $V$, we write $\textrm{Coh}(\mathcal O_V)$ for the category of coherent $\mathcal O_V$-modules on $V$.

\begin{cor} \label {coherent}
Let $(V, \mathcal O_V)$ be a locally compact ringed space that is countable at infinity with $\mathcal O_V$ coherent.
Let $i : T \hookrightarrow V$ the inclusion of a closed subspace.
Then, the restriction functors $\textrm{Coh}(\mathcal O_{V'}) \to \textrm{Coh}(i^{-1}\mathcal O_V)$ where $V'$ runs trough all inclusions of open neighborhoods of $T$ in $V$ induce an equivalence of categories
$$
\varinjlim \textrm{Coh}(\mathcal O_{V'}) \simeq \textrm{Coh}(i^{-1}\mathcal O_V).
$$
\end{cor}

\textbf {Proof : }Using the previous proposition, this is an immediate consequence of Proposition \ref{limit2}. $\Box$\bigskip

We will now apply these results to our analytic varieties.
We have the first fundamental result:

\begin{prop}
If $(X, V)$ is an analytic variety over $\mathcal V$, then $i_X^{-1}\mathcal O_V$ is a coherent ring.
\end{prop}

{\bf Proof :} This is a direct application of Proposition \ref{icoher}. $\Box$\bigskip

The next step consists in applying the results of corollary \ref{coherent}.
We will say that an analytic variety $(X, V)$ is \emph{countable at infinity} if there exists a strict neighborhood $V'$ of $X$ in $V$ which is countable at infinity with $]X[_{V}$ closed in $V'$.

\begin{prop}

\begin{enumerate} \renewcommand{\labelenumi}{\roman{enumi})}

\item If $X \subset P$ is a formal embedding with $P$ quasi-compact, then $(X, P_K)$ is countable at infinity.

\item Locally, any good analytic variety $(X, V)$ is countable at infinity.
\end{enumerate}

\end{prop}

\textbf {Proof : } We prove the first assertion : we have $]X[_P$ closed in $]\bar X[_P$ and $]\bar X[_P$ open in $P_K$. It is therefore sufficient to recall that the tube $]\bar X[_P$ is the increasing union of the closed tubes $[\bar X]_{P\eta}$, that these closed tubes are closed subsets of $P_K$, which is compact, and are therefore compact themselves.

Now, we prove the second assertion : we may first shrink $V$ in order to have $]X[_V$ closed in $V$, and then use the fact that any good variety has a basis of open subsets countable at infinity. $\Box$\bigskip

\begin{prop} \label {coherent2}
If $(X, V)$ is an analytic variety over $S$ which is countable at infinity, we have an equivalence of categories
$$
\varinjlim \textrm{Coh}(\mathcal O_{V'}) \simeq \textrm{Coh}(i_X^{-1}\mathcal O_{V})
$$
when $V'$ runs trough all inclusions of open neighborhoods of $]X[_V$ in $V$.
\end{prop}

\textbf {Proof : }This is a particular case of corollary \ref {coherent}. $\Box$\bigskip

%
%
\section{Crystals}

With the use of Berkovich theory, we do not need Berthelot's $j^\dagger$ construction as the next result shows.

\begin{prop} \label{jdagger}
Let $(X, V)$ be an analytic variety over $\mathcal V$ and $\mathcal F$ any sheaf on $V$.
If $]X[_{V}$ is closed in $V$, we have
$$
j^\dagger \mathcal F := \varinjlim j'{_*}j^{'-1} \mathcal F = i_*i^{-1} \mathcal F,
$$
where $j'$ runs through all strict neighborhoods of $X$ in $V$.
\end{prop}

\textbf {Proof : } 
It is true in general that if $V$ is a topological space whose points are closed, then any subset $T$ of $V$ is an intersection of open subsets.
It follows that if $T$ is closed in $V$ and $\mathcal F$ any sheaf, then
$$
j^\dagger \mathcal F := \varinjlim j'{_*}j^{'-1} \mathcal F = i_*i^{-1} \mathcal F,
$$
where $j' : V' \hookrightarrow V$ runs through all inclusions of open neighborhoods of $T$ in $V$ and $i : T \hookrightarrow V$ denotes the inclusion map.
This is easily checked by looking at the stalks. $\Box$\bigskip

\begin{prop} \label{dag=-1}
\begin{enumerate} \renewcommand{\labelenumi}{\roman{enumi})}

\item Any morphism $(f, u) : (X', V') \to (X, V)$ of analytic varieties over $\mathcal V$ induces a morphism of ringed spaces
$$
(]f[_{u}^\dagger, ]f[_{u*}) : (]X'[_{V'}, i_{X',V'}^{-1} \mathcal O_{V'}) \to (]X[_V, i_{X,V}^{-1} \mathcal O_V).
$$
and we have for all $i_{X,V}^{-1} \mathcal O_V$-module $\mathcal F$,
$$
]f[_{u}^\dagger \mathcal F = i_{X',V'}^{-1} u^* i_{X,V*} \mathcal F.
$$
\item Let $(X, V)$ be an analytic variety, $\mathcal F$ a $i_{X,V}^{-1}\mathcal O_{V}$-module and $\alpha : Y \hookrightarrow X$ the inclusion of a subscheme, then $]\alpha[_V^\dagger\mathcal F = ]\alpha[_V^{-1}\mathcal F$.
 
\end{enumerate}

\end{prop}

\textbf {Proof : } Pulling back by $i_{X',V'}$ the canonical map $u^{-1}\mathcal O_V \to \mathcal O_{V'}$ gives
$$
]f[_u^{-1}i_{X,V}^{-1}\mathcal O_V = i_{X',V'}^{-1}u^{-1}\mathcal O_V \to i_{X',V'}^{-1}\mathcal O_{V'}
$$
and we get a morphism of ringed spaces as asserted.
If $\mathcal F$ is a sheaf on $]X[_V$, we have 
$$
]f[_u^{-1}\mathcal F = ]f[_u^{-1} i_{X,V}^{-1}i_{X,V*}\mathcal F = i_{X',V'}^{-1}u^{-1}i_{X,V*}\mathcal F.
$$
Thus, if $\mathcal F$ is a $i_{X,V}^{-1} \mathcal O_V$-module, we see that
$$
]f[_{u}^\dagger \mathcal F = i_{X',V'}^{-1}\mathcal O_{V'} \otimes_{i_{X',V'}^{-1}u^{-1}\mathcal O_V} i_{X',V'}^{-1}u^{-1}i_{X,V*}\mathcal F
$$
$$
= i_{X',V'}^{-1}[\mathcal O_{V'} \otimes_{u^{-1}\mathcal O_V} u^{-1}i_{X,V*}\mathcal F] = i_{X',V'}^{-1} u^* i_{X,V*} \mathcal F.
$$
The formula in the second assertion comes from
$$
]\alpha[_V^\dagger \mathcal F = i_{X',V}^{-1} i_{X,V*} \mathcal F = ]\alpha[_V^{-1} i_{X,V}^{-1} i_{X,V*} \mathcal F = ]\alpha[_V^{-1} \mathcal F
$$ $\Box$\bigskip

Actually, in practice, we will often write $u^\dagger$ and $u_{*}$ instead of $]f[_{u}^\dagger$ and $]f[_{u*}$.
On the other hand, when $u = \textrm{Id}_V$, we will write $]f[_V$ as we actually did in the second part of the proposition.
We will also write $ i_X^{-1} \mathcal O_V := i_{X,V}^{-1} \mathcal O_V$.

\begin{cor}
The presheaf of rings
$$
\mathcal O_{\mathcal V}^\dagger : (X, V) \mapsto \Gamma(]X[_V, i_X^{-1}\mathcal O_V)
$$
is a sheaf on $\textrm{AN}^\dagger(\mathcal V)$.
\end{cor}

\textbf {Proof : }
Follows from the first part of the proposition. $\Box$\bigskip

This is the \emph{sheaf of overconvergent functions on $\mathcal V$}.
If $T$ is an overconvergent presheaf on $\mathcal V$, the localization $\mathcal O_{T}^\dagger$ of $\mathcal O_{\mathcal V}^\dagger$ to $\textrm{AN}^\dagger(T)$ is the \emph{sheaf of overconvergent functions on $T$}. A module over this ring will be called an \emph{overconvergent module} on $T$.
By definition, the realization of $\mathcal O_{\mathcal V}^\dagger$ on some analytic variety $(X, V)$ is just $i_{X,V}^{-1} \mathcal O_V$.
Note also that if $f : T' \to T$ is a morphism of analytic presheaves on $\mathcal V$, we have
$$
f_{\textrm{AN}^\dagger}^{-1}\mathcal O_{T}^\dagger = \mathcal O_{T'}^\dagger.
$$
If $E$ is an overconvergent module on $(X, V)$, then
$E_V$ is a $ i_X^{-1} \mathcal O_V$-module and for all morphisms
$$
(f, u) : (X', V') \to (X, V),
$$
the morphism $\phi_{uE} : ]f[_u^{-1}E_V \to E_{V'}$ on $]X'[_{V'}$ extends to a $i_X^{-1} \mathcal O_{V'}$-linear map
$$
\phi^\dagger_{fuE} : ]f[_u^\dagger E_V \to E_{V'}.
$$
\begin{prop} \label{adjunction}
Let $(X, V)$ be an analytic variety over $\mathcal V$. Then,
\begin{enumerate}

\item There is a canonical morphism of ringed spaces
$$
(\varphi_{X,V}^\dagger, \varphi_{X,V*}) : (\textrm{AN}^\dagger(X,V), \mathcal O_{(X,V)}^\dagger) \to (]X[_V, i_X^{-1}\mathcal O_V).
$$
\item If $\mathcal F$ is any $i_X^{-1}\mathcal O_V$-module and $(f, u) : (X', V') \to (X, V)$ is any morphism, then the realization of $\varphi_{X,V}^\dagger\mathcal F$ on $(X',V')$ is just $]f[_u^\dagger \mathcal F$.

\item If $E$ is an overconvergent module on $(X,V)$, then the realization of the adjunction map
$$
\varphi_{X,V}^\dagger \varphi_{X,V*}E \to E
$$
along some $(f, u) : (X', V') \to (X, V)$ is the transition map $\phi^\dagger_{uE} : ]f[_u^\dagger E_V \to E_{V'}$.

\end{enumerate}

\end{prop}

\textbf {Proof : }We know that
$$
i_X^{-1}\mathcal O_V = \varphi_{X,V*}\mathcal O_{(X,V)}^\dagger.
$$
The first assertion follows formally from this fact.
Now, by definition, we have
$$
\varphi_{X,V}^\dagger\mathcal F = \mathcal O_{(X,V)}^\dagger \otimes_{\varphi_{X,V}^{-1} i_X^{-1}\mathcal O_V} \varphi_{X,V}^{-1}\mathcal F.
$$
Since realization on $(X', V')$ is the the pull-back by $\psi_{X',V'}$, it commutes to tensor products and we see that
$$
(\varphi_{X,V}^\dagger\mathcal F)_{X',V'} = i_{X'}^{-1}\mathcal O_{V'} \otimes_{]f[_u^{-1} i_X^{-1}\mathcal O_V} ]f[_u^{-1}\mathcal F = ]f[_u^{\dagger}\mathcal F.
$$
The case $\mathcal F = E_V$ gives the last result. $\Box$\bigskip

\begin{prop}
\label{carac} If $T$ is an overconvergent presheaf on $\mathcal V$, the category of overconvergent modules on $T$ is equivalent to the following category :
\begin{enumerate}
\item An object is a collection of $i_{X}^{-1} \mathcal O_{V}$-modules $E_{X,V}$ on $]X[_{V}$ for each $(X, V) \in \textrm{AN}^\dagger(T)$ and $i_{X'}^{-1} \mathcal O_{V'}$-linear maps $\phi^\dagger_{fuE} : ]f[_u^\dagger E_{X,V} \to E_{X',V'}$ for each morphism
$$
(f, u) : (X', V') \to (X, V) 
$$
of analytic varieties over $T$, satisfying the usual compatibility conditions.
\item A morphism is a collection of $i_{X}^{-1} \mathcal O_{V}$-linear maps $E_{X,V} \to E'_{X,V}$ compatible with the morphisms $\phi^\dagger_{fu}$.
\end{enumerate}
\end{prop}

\textbf {Proof : } As usual. $\Box$\bigskip

Let $T$ be an overconvergent presheaf on $\mathcal V$.
An overconvergent module $E$ on $T$ is an \emph{overconvergent crystal} if all the transition maps $\phi^\dagger_{fuE}$ in proposition \ref{carac} are isomorphisms.
We denote this full subcategory by $\textrm{Cris}^\dagger(T/S)$.
Note that an overconvergent module $E$ on $T$ is a crystal if and only if for all $(X,V) \in \textrm{AN}^\dagger(T)$, $E_{/(X,V)}$ is an overconvergent crystal on $(X,V)$.
It is actually sufficient to check it for an analytic covering $\{(X_i, V_i)\}_{i \in I}$ of $T$.
We should also remark that any morphism of analytic presheaves $f : T' \to T$ provides a functor
$$
f^{-1}_{\textrm{AN}} : \textrm{Cris}^\dagger(T') \to \textrm{Cris}^\dagger(T).
$$
\begin{prop} \label {nice}
If $(X, V)$ is an analytic variety over $\mathcal V$, the functors
$\varphi_{X,V}^\dagger$ and $\varphi_{X,V*}$ induce an equivalence of categories between $\textrm{Cris}^\dagger(X,V)$ and the category of $i_X^{-1}\mathcal O_V$-modules on $]X[_V$.
In particular, $\textrm{Cris}^\dagger(X,V)$ is an abelian category with tensor product, internal hom and enough injectives.
\end{prop}

\textbf {Proof : }It follows from the second assertion of \ref{adjunction}\spaceÊthat if $\mathcal F$ is any $i_X^{-1}\mathcal O_V$-module, then $\varphi_{X,V}^\dagger \mathcal F$ is a crystal and the adjunction map $\varphi_{X,V*}\varphi_{X,V}^\dagger \mathcal F \to \mathcal F$ is bijective.
Now, if $E$ is an overconvergent crystal on $(X, V)$, we know from the last assertion of \ref{adjunction}\spaceÊthat the realization of the adjunction map
$\varphi_{X,V}^\dagger \varphi_{X,V*}E \to E$ along any $(f, u) : (X', V') \to (X, V)$ is the transition map $\phi^\dagger_{uE} : ]f[_u^\dagger E_V \to E_{V'}$, which is by hypothesis, an isomorphism.
It follows that this adjunction map is also an isomorphism. $\Box$\bigskip

\begin{cor} \label{catT}
If $T$ is an overconvergent presheaf on $\mathcal V$, the category $\textrm{Cris}^\dagger(T)$ is an additive subcategory of the category of overconvergent modules which is stable under cokernel, extensions and tensor product.
\end{cor}

\textbf {Proof : } We may clearly assume that $T = (X,V)$ in which case everything follows from the right exactness of $\varphi^\dagger_{X,V}$. $\Box$\bigskip

\begin{prop}
Let $T$ be an overconvergent presheaf on $\mathcal V$.
If $E$ is an overconvergent crystal on $T$ and $E'$ is any overconvergent module, then for each analytic variety $(X,V)$ over $T$, we have
$$
\mathcal H \textrm{om}_{\mathcal O_{T}^\dagger}(E,E')_{X,V} =\mathcal H\textrm{om}_{i_{X}^{-1}\mathcal O_{V}}(E_{X,V}, E'_{X,V}).
$$
\end{prop}

\textbf {Proof : }It is again sufficient to consider the case $T = (X,V)$.
Then, our assertion formally follows from Proposition \ref{nice}. Namely, we have 
$$
\mathcal H\textrm{om}_{\mathcal O_{S/X,V}^\dagger}(E,E')_{X,V} = \varphi_{X,V*}\mathcal H\textrm{om}_{\mathcal O_{S/X,V}^\dagger}(\varphi_{X,V}^\dagger \varphi_{X,V*}E, E')
$$
$$
= \mathcal H\textrm{om}_{\varphi_{X,V*}\mathcal O_{S/X,V}^\dagger}(\varphi_{X,V*}E, \varphi_{X,V*}E') = \mathcal H\textrm{om}_{i_{X}^{-1}\mathcal O_{V}}(E_{X,V}, E'_{X,V}). \quad \Box
$$\bigskip

At this point, we need to introduce finiteness conditions.

\begin{prop} \label {finite}
Let $T$ be an overconvergent presheaf over $\mathcal V$. An overconvergent module $E$ on $T$ is finitely presented if and only if it is a crystal and for all analytic varieties $(X,V)$ over $T$, $E_V$ is a coherent $i_X^{-1}\mathcal O_V$-module.
Moreover, $E$ is locally free of rank $r$ if and only if for all $(X,V)$ over $T$, $E_V$ is a locally free $i_X^{-1}\mathcal O_V$-module of rank r.
\end{prop}

\textbf {Proof : } The question being local on $\textrm{AN}(T)$, we may assume that $T = (X,V)$.
Then, since both $\varphi_{X,V*}$ and $\varphi_{X,V}^\dagger$ are additive and right exact, our assertions follow directly from the fact that $\mathcal O^\dagger_{X,V}$ is a crystal. $\Box$\bigskip

If $T$ is an overconvergent presheaf on $\mathcal V$, we will denote by $\textrm{Mod}_{\textrm{fp}}^\dagger(T)$ the category of overconvergent modules of finite presentation on $T$.

\begin{prop}
Let $T$ be an overconvergent presheaf on $\mathcal V$.
If $E, E' \in \textrm{Mod}_{\textrm{fp}}^\dagger(T)$, then
$$
\mathcal H\textrm{om}_{\mathcal O_{T/S}^\dagger}(E,E') \in \textrm{Mod}_{\textrm{fp}}^\dagger(T).
$$
In particular, it is a crystal.
\end{prop}

\textbf {Proof : }If $(f, u) : (X', V') \to (X, V)$ is a morphism over $T$, we have
$$
]f[_u^\dagger\mathcal H\textrm{om}_{\mathcal O_{T}^\dagger}(E,E')_{X,V} = ]f[_u^\dagger\mathcal H\textrm{om}_{i_{X}^{-1}\mathcal O_{V}}(E_{X,V}, E'_{X,V})
$$
and since $E_{X,V}$ is finitely presented, we have
$$
]f[_u^\dagger\mathcal H\textrm{om}_{i_{X}^{-1}\mathcal O_{V}}(E_{X,V}, E'_{X,V})
= \mathcal H\textrm{om}_{i_{X'}^{-1}\mathcal O_{V}}(]f[_u^\dagger E_{X,V}, ]f[_u^\dagger E'_{X,V})
$$
and since we are dealing with crystals, we have 
$$
\mathcal H\textrm{om}_{i_{X'}^{-1}\mathcal O_{V}}(]f[_u^\dagger E_{X,V}, ]f[_u^\dagger E'_{X,V}) = \mathcal H\textrm{om}_{i_{X'}^{-1}\mathcal O_{V}}(E_{X',V'}, E'_{X',V'})
$$
and we get
$$
]f[_u^\dagger\mathcal H\textrm{om}_{\mathcal O_{T/S}^\dagger}(E,E')_{X,V} = \mathcal H\textrm{om}_{\mathcal O_{T}^\dagger}(E,E')_{X',V'}. \quad \Box
$$

We will finish this section with the study of immersions of algebraic varieties.
We need some preliminary results.

\begin{lem}
Let $i : T \hookrightarrow V$ be the inclusion of a subspace in a topological space, $\mathcal A$ a sheaf of rings on $V$, and $\mathcal F$ and $\mathcal G$ two $\mathcal A$-modules. Then 
$$
i_*i^{-1}(\mathcal F \otimes_\mathcal A \mathcal G) = i_*i^{-1}\mathcal F \otimes_\mathcal A \mathcal G.
$$
\end{lem}

\textbf {Proof : } Since tensor product commutes to pull-back, this immediately follows from the fact that $i_*$ is fully faithful. $\Box$\bigskip

The following lemma is analogous to Proposition 2.1.4 of \cite {Berthelot96}.

\begin{lem} \label{commute}
Let $(f, u) : (X', V') \to (X, V)$ be a morphism of analytic varieties over $\mathcal V$.
Assume that $]X[_V$ is closed in $V$ and that $u_K^{-1}(]X[_V) = ]X'[_{V'}$.
Let us write $i : ]X[_V \hookrightarrow V$ and $i' : ]X'[_{V'} \hookrightarrow V'$ for the inclusion maps. Then, if $\mathcal F$ is any $\mathcal O_V$-module, we have
$$
u^*i_*i^{-1} \mathcal F \simeq i'_*i'^{-1}u^*\mathcal F.
$$
In particular, if $\mathcal F$ is a $i_X^{-1} \mathcal O_V$-module, we have
$$
u^*i_* \mathcal F \simeq i'_*]f[_u^\dagger \mathcal F.
$$
\end{lem}

\textbf {Proof : }Using the previous lemma, our statement is a formal consequence of the fact that, with our hypothesis (a cartesian square with a closed embedding), we have for any sheaf, $u^{-1}i_{*} \mathcal F = i'_{*}u^{-1} \mathcal F$. $\Box$\bigskip

\begin{prop} Let $S$ be a formal $\mathcal V$-scheme and $\alpha : Y \hookrightarrow X$ an open immersion of $S_k$-schemes. If $E$ is an overconvergent crystal on $Y/S$ then $\alpha_{\textrm{AN}*}E$ is an overconvergent crystal on $X/S$.
In other words, $\alpha_{\textrm{AN}*}$ induces a functor
$$
\alpha_{\textrm{AN}*} : \textrm{Cris}^\dagger(Y/S) \to \textrm{Cris}^\dagger(X/S).
$$
\end{prop}

\textbf {Proof : }We are given an overconvergent crystal $E$ on $Y/S$, and a morphism $(f, u) = (X' ,V') \to (U, V)$ in $\textrm{AN}(X/S)$.
We want to show that
$$
]f[_u^\dagger( \alpha_{\textrm{AN}^*}E)_V = ( \alpha_{\textrm{AN}*}E)_{V'}.
$$
We may clearly assume that $U = X$. We call $Y'$ the pull-back of $Y$ in $X'$, and denote by $\alpha' : Y' \hookrightarrow X'$ the inclusion map.
Using Proposition \ref {immer}, our condition can be rewritten
$$
]f[_u^\dagger ]\alpha[_{V*}E_{Y,V} = ]\alpha'[_{V'*} E_{Y',V'}.
$$
But if $g : Y' \to Y$ denotes the map induced by $f$, we have $E_{Y',V'} = ]g[_u^\dagger E_{Y,V}$.
We are thus reduced to show that
$$
]f[_u^\dagger \circ ]\alpha[_{V*} = ]\alpha'[_{V'*} \circ ]g[_u^\dagger.
$$
We may split the verification in two parts. So we assume first that $u = \textrm{Id}_V$. We simply have here $]f[^\dagger = ]f[^{-1}$ and $]g[^\dagger = ]g[^{-1}$. Thus, we are reduced to check that $]f[^{-1} \circ ]\alpha'[_{*} = ]\alpha[_{*} \circ ]g[^{-1}$ which follows from the fact that $]Y]_V$ is closed in $]X[_V$ (and the diagram is cartesian).

We assume now that $P' = P$ and $f = \textrm{Id}_X$, in which case, also $g = \textrm{Id}_Y$. Shrinking $V$ and $V'$ if necessary, we may also assume that $]X[_V$ and $]X[_{V'}$ are closed in $V$.
Moreover, since $i_{X,V'*}$ is fully faithful, it is sufficient to prove that
$$
i_{X,V'*} \circ ]\textrm{Id}_X[_u^\dagger \circ ]i[_{V*} = i_{X,V'*} \circ ]\alpha'[_{V'*} \circ ]\textrm{Id}_Y[_u^\dagger.
$$
We are in the situation of applying lemma \ref{commute} both to $X$ and $Y$.
On the left hand side, we get
$$
i_{X,V'*} \circ ]\textrm{Id}_X[_u^\dagger \circ ]\alpha[_{V*} = u^* \circ  i_{X,V*} \circ ]\alpha[_{V*} = u^* \circ  i_{Y,V*}
$$
and on the 
On the right hand side, we get
$$
 i_{X,V'*} \circ ]\alpha'[_{V'*} \circ ]\textrm{Id}_Y[_u^\dagger = i_{Y,V'*} \circ ]\textrm{Id}_Y[_u^\dagger = u^* \circ  i_{Y,V*}. \quad \Box
$$
\bigskip

%
%
\section{Stratifications}

In this section (and the next one), we fix a formal $\mathcal V$-scheme $S$. We will use the notion of stratification as a bridge between crystals and modules with integrable connection.

\bigskip

Let $(X, V)$ an analytic variety over $S$. We embed $X$ diagonally in $P^2 := P \times_{S} P$ and denote by
$$
p_1, p_2 : V^2 := V \times_{S_K} V \to V
$$
the projections. 
An \emph{overconvergent stratification }Êon a $i_X^{-1}\mathcal O_V$-module $\mathcal F$ is an isomorphism
$$
\epsilon : p_2^\dagger \mathcal F \simeq p_1^\dagger \mathcal F
$$
on $]X[_{V^2}$, called the \emph{Taylor isomorphism} of $E$, satisfying the usual cocycle condition on triple products.
A \emph{morphism of overconvergent stratified modules} is a morphism of $i_X^{-1}\mathcal O_V$-modules compatible with the data.
We denote this category by $\textrm{Strat}^\dagger(X, V/S)$.
We will also set :
$$
\mathcal H^{\dagger}(\mathcal F) := \ker \left[\epsilon \circ p_2^{-1} - p_1^{-1} : p_{]X[_V/S*} \mathcal F \to p_{]X[_{V^{2}}/S*} p_1^\dagger \mathcal F\right].
$$
\begin{lem}
If $(f, u) : (X', V') \to (X, V)$ is a morphism of analytic varieties over $\mathcal V$ with $u$ flat in a neighborhood of $]X'[_{V'}$, then $]f[_u^\dagger$ is exact.
\end{lem}

\textbf {Proof : } Our assumtions is that that $u^{-1}\mathcal O_V \to \mathcal O_{V'}$ is flat.
Pulling back by $ i_{X'}^{-1}$, we see that
$$
]f[_u^{-1}i_X^{-1}\mathcal O_V = i_{X'}^{-1} u^{-1} \mathcal O_V \to i_{X'}^{-1}\mathcal O_{V'}
$$
is also flat. And this is what we want. $\Box$\bigskip

\begin{prop} \label {flatness}
Let $(X, V)$ an analytic variety over $S$.
Then, the category $\textrm{Strat}^\dagger(X, V/S)$ is an additive category with cokernels and the forgetful functor to from $\textrm{Strat}^\dagger(X, V/S)$ to the category of $i_X^{-1}\mathcal O_V$-modules is right exact and faithful. Moreover, if $V$ is universally flat in a neighborhood of $]X[_V$ over $S_K$, then $\textrm{Strat}^\dagger(X, V/S)$ is even an abelian category and the forgetful functor is exact.
\end{prop}

\textbf {Proof : }It is clear that we have an additive category and that the forgetful functor is faithful. We will show, when $V$ is universally flat, the existence of a stratification on the kernel of a morphism of dagger stratified modules $m : \mathcal F \to \mathcal G$, and also that this new structure turns $\ker m$ into a kernel in the category of dagger stratified modules. The analogous result for cokernels is showed exactly in the same way (without the flatness assumption).

If $V$ is universally flat over $S$, the projections $p_1$ and $p_2$ are flat. It follows from the lemma that $p_1^\dagger$ and $p_2^\dagger$ are exact. Thus, we see that $\epsilon_\mathcal F$ induces an isomorphism $p_2^\dagger(\ker m) \simeq p_1^\dagger(\ker m)$ which is obviously a stratification on $\ker m$. Clearly, if the image of a morphism of dagger stratified modules $\mathcal H \to \mathcal F$ is contained in $\ker m$, the induced map $\mathcal H \to \ker m$ is compatible with the stratifications. $\Box$\bigskip

We will need an infinitesimal version of this notion of dagger stratification. We want first to recall some notions concerning analytic varieties. If $W$ is a fixed analytic space over $K$ and $V$ an analytic space over $W$, we will write $V^{(n)}$ for the $n$-th infinitesimal neighborhood of $V$ in $V \times_{W} V$ and $p_1^{(n)}, p_2^{(n)} : V^{(n)} \to V$ for the projections. By definition, there is an exact sequence
$$
0 \to \Omega^1_{V/W} \to \mathcal O_{V^{(2)}} \to \mathcal O_V \to 0.
$$
We can define stratifications, modules with (integrable) connexion and the sheaf of differential operators $\mathcal D_{V/W}$ as usual. Everything behaves as expected.

Now, if $(X \subset P \leftarrow V) \in \textrm{AN}^\dagger(S)$, we may take $W = S_K$ and consider
$$
(X \subset P \times_S P \leftarrow V^{(n)}) \in \textrm{AN}(S).
$$
A \emph{stratification} on a $i_X^{-1}\mathcal O_V$-module $\mathcal F$ is a compatible sequence of \emph{Taylor isomorphisms}
$$
\{\epsilon^{(n)} : p_2^{(n)\dagger} \mathcal F \simeq p_1^{(n)\dagger} \mathcal F\}_{n \in \mathbf N}
$$
that satisfy the cocycle condition on triple products (and $\epsilon^{(0)} = Id_\mathcal F$). A morphism of stratified modules is a morphism of $i_X^{-1}\mathcal O_V$-modules that is compatible with the data.
We will write 
$$
\mathcal H^{(n)}(\mathcal F) := \ker \left[(\epsilon^{(n)} \circ (p_2^{(n)})^{-1}) - (p_1^{(n)})^{-1} : p_{]X[_V/S*} \mathcal F \to p_{]X[_{V}/S*} p_1^{(n)\dagger} \mathcal F\right].
$$
A \emph{connexion} on a $i_X^{-1}\mathcal O_V$-module $\mathcal F$ is an $\mathcal O_{S_K}$-linear map
$$
\nabla : \mathcal F \to \mathcal F \otimes_{i_X^{-1}\mathcal O_V} i_X^{-1}\Omega^1_{V/S_K}
$$
satisfying the Leibniz rule.
A \emph{horizontal} map is a $i_X^{-1}\mathcal O_V$-linear map compatible with the connexions. \emph{Integrability} is defined in the usual way.

\bigskip

Stratified modules $i_X^{-1}\mathcal O_V$-modules form a category $\textrm{Strat}(X, V/S)$ and $i_X^{-1}\mathcal O_V$-modules with integrable connexion make a category $\textrm{MIC}(X, V/S)$.
As usual, there is an obvious forgetful functor $\textrm{Strat}(X, V/S) \to \textrm{MIC}(X, V/S)$ sending $\{\epsilon^{(n)}\}_{n \in \mathbf N}$ to the morphism
$$
\nabla = (\epsilon^{(1)} \circ (p_2^{(1)})^{-1}) - (p_1^{(1)})^{-1} : \mathcal F \to \mathcal F \otimes_{i_X^{-1}\mathcal O_V} i_X^{-1}\Omega^1_{V/S_K} \subset p_1^{(2)\dagger} \mathcal F.
$$
When $\textrm{Char} K = 0$ and $V$ is smooth in the neighborhood of $]X[_V$, it is not hard to see that there is an equivalence
$$
\textrm{Strat}(X, V/S) \simeq \textrm{MIC}(X, V/S)
$$
and that $\mathcal H^{(n)}(\mathcal F)$ is independent of $n$, for $n > 0$, and canonically isomorphic to $\mathcal F^{\nabla = 0}$.

\bigskip

We can do better with some finiteness conditions.
In general, if $\mathcal F$ and $\mathcal G$ are two $i_{X}^{-1}\mathcal O_{V}$-modules with integrable connexions, then $\mathcal H\textrm{om}_{i_{X}^{-1}\mathcal O_{V}}(\mathcal F, \mathcal G)$ has a natural integrable connexion given by $\nabla (m)(x) = \nabla(m(x)) - (m \otimes \textrm{Id})(x)$. Moreover, we have
$$
\Gamma(]X[_V, \mathcal H\textrm{om}_{i_{X}^{-1}\mathcal O_{V}}(\mathcal F, \mathcal G)^{\nabla = 0} ) = \textrm{Hom}_{\textrm{MIC}(X,V/S)}(\mathcal F, \mathcal G).
$$
In the same way, if $\mathcal F$ and $\mathcal G$ be coherent stratified $i_{X}^{-1}\mathcal O_{V}$-modules, then $\mathcal H\textrm{om}_{i_{X}^{-1}\mathcal O_{V}}(\mathcal F, \mathcal G)$ has a natural stratification and we have
$$
\varprojlim \Gamma(]X[_V, \mathcal H^{(n)}[\mathcal H\textrm{om}_{i_{X}^{-1}\mathcal O_{V}}(\mathcal F, \mathcal G)]) = \textrm{Hom}_{\textrm{Strat}(X,V/S)}(\mathcal F, \mathcal G)
$$
Moreover, the canonical functor
$$
\textrm{Strat}(X, V/S) \to \textrm{MIC}(X, V/S)
$$
is compatible with these constructions.
We now come back to overconvergent stratifications:
By restriction, there is an obvious functor
$$
\textrm{Strat}^\dagger(X, V/S) \to \textrm{Strat}(X, V/S)
$$
and a compatible family of canonical maps $\mathcal H^\dagger(\mathcal F) \to \mathcal H^n(\mathcal F)$.
Assume that $\mathcal F$ and $\mathcal G$ are two coherent dagger stratified $i_{X}^{-1}\mathcal O_{V}$-modules. Then $\mathcal H\textrm{om}_{i_{X}^{-1}\mathcal O_{V}}(\mathcal F, \mathcal G)$ has a natural dagger stratification and
$$
\Gamma(]X[_V, \mathcal H^{\dagger}[\mathcal H\textrm{om}_{i_{X}^{-1}\mathcal O_{V}}(\mathcal F, \mathcal G)]) = \textrm{Hom}_{\textrm{Strat}^\dagger(X,V/S)}(\mathcal F, \mathcal G)
$$
Again, the canonical functor
$$\
\textrm{Strat}^\dagger(X, V/S) \to \textrm{Strat}(X, V/S)
$$
is compatible with these constructions.

\bigskip

In the next proposition, we will need to localize with respect to the Grothendieck topology of an analytic variety $V$.
It might be convenient to denote by $V_{G}$ the corresponding Grothendieck space and by $\pi_{V} : V_{G} \to V$ the canonical morphism of ringed spaces.
Note that, by functoriality, if $i : T \hookrightarrow V$ is the inclusion of an analytic domain and $\mathcal F$ a sheaf on $V$, then
$$
\Gamma(T, i^{-1}\mathcal F) = \Gamma(T, \pi_{V}^{-1}\mathcal F).
$$
If $\mathcal F$ is an $\mathcal O_{V}$-module, it is common to write $\mathcal F_{G} = \pi_{V}^*\mathcal F$.
When $V$ is good, we get an equivalence on coherent sheaves.
\bigskip

\begin{prop} \label{goodfaith}
If $V$ is a good analytic variety and $\mathcal F$ is a coherent sheaf on $V$, then the canonical map $\pi_{V}^{-1}\mathcal F \to \mathcal F_{G}$ is injective.
\end{prop}

\textbf {Proof : }
It is sufficient to show that if $W$ is an affinoid domain inside $V$, the canonical map
$$
\Gamma(W, \pi_{V}^{-1}\mathcal F) \to \Gamma(W, \mathcal F_{G})
$$
is injective.
By definition, we have
$$
\Gamma(W, \pi_{V}^{-1}\mathcal F) = \varinjlim_{W \subset U}\Gamma(U, \mathcal F)
$$
where $U$ runs through analytic open subsets of $V$.
It is therefore sufficient to prove that if $U$ is an open neighborhood of $W$ in $V$ and $s \in \Gamma(U, \mathcal F)$ is sent to $0 \in \Gamma(W, \mathcal F)$, then there exists an open neighborhood $U'$ of $W$ in $V$ such that $s_{|U'} = 0$.
Of course, this will follow if we show that for each $x \in W$, the stalk of $s$ at $x$ is zero.
But if we denote by $\mathcal F_{W}$ the coherent sheaf on $W$ induced by $\mathcal F$ (i.e. $\mathcal F_{W} := \pi_{W*}(\mathcal F_{G})_{|W}$), the image of $s$ in the stalk of $\mathcal F_{W}$ at $x$ is zero.
It is therefore sufficient to check that the map
$$
\mathcal F_{x} \to \mathcal F_{Wx}
$$
is injective.
Since $V$ is good, we may assume that $V$ and $W$ are both affinoid with $A := \Gamma(V, \mathcal O_{V})$ and $M := \Gamma(V, \mathcal F)$.
And we want to show that the map
$$
\mathcal O_{V,x} \otimes_{A} M \to \mathcal O_{W,x} \otimes_{A} M
$$
is injective.
This is an immediate consequence of the flatness of the inclusion map $W \hookrightarrow V$. $\Box$

\begin{prop} Let $X \hookrightarrow P$ be a good admissible formal embedding into a quasi-compact formal $S$-scheme which is smooth at $X$.
If $\mathcal F$ is a coherent overconvergent stratified module on $(X \subset P)/S$, then
$$
\mathcal H^{\dagger}(\mathcal F) =\varprojlim \mathcal H^{(n)}(\mathcal F).
$$
\end{prop}

\textbf {Proof : }
In each case, we denote by $p_{1} : ]X[_{P^2} \to ]X[_{P}$ and $p_{1}^{(n)} : ]X[_{P^{(n)}} \to ]X[_{P}$ the first projection.
By definition, it is sufficient to show that the natural map
$$
p_{]X[_{P^2}/S*} p_1^{\dagger} \mathcal F \to \varprojlim p_{]X[_{P}/S*} p_1^{(n)\dagger} \mathcal F
$$
is injective.
It is even sufficient to check that the canonical map
$$
p_{1*} p_1^{\dagger} \mathcal F \to \varprojlim p_1^{(n)\dagger} \mathcal F
$$
is injective.
Since $\mathcal F$ is coherent and $(X, P_{K})$ countable at infinity, it is sufficient to consider sheaves of the form $i^{-1} \mathcal F$ where $\mathcal F$ is a coherent module on some good smooth strict neighborhood of $X$ in $V$ and $i : ]X[_{P} \hookrightarrow V$ denotes the inclusion map.
In this situation, we want to show that
$$
p_{1*} {i'}^{-1} p_1^{*} \mathcal F \to \varprojlim i^{-1} p_1^{(n)*} \mathcal F
$$
is injective with $i' : ]X[_{P^2} \hookrightarrow V^2$.

Now, we want to use the Grothendieck topology but we need to be very careful because $]X[_{P}$ might not be good.
Anyway, since $V$ is good, $\mathcal F$ extends uniquely to a coherent sheaf $\mathcal F_{G} = \pi_{V}^*\mathcal F$ for the Grothendieck topology of $V$.
If $W$ is an analytic open subset of $]X[_{P}$, since $V^2$ is good and $p_1^{*} \mathcal F_{G}$ coherent, it follows from lemma \ref{goodfaith} that there is a natural injective map
$$
\Gamma(W, p_{1*} {i'}^{-1} p_1^{*} \mathcal F) = \Gamma(p_{1}^{-1}(W), {i'}^{-1} p_1^{*} \mathcal F)
$$
$$
= \Gamma(p_{1}^{-1}(W), \pi_{V^2}^{-1} p_1^{*} \mathcal F) \hookrightarrow \Gamma(p_{1}^{-1}(W), \pi_{V^2}^{*} p_1^{*} \mathcal F)
$$
$$
= \Gamma(W, p_{1*} p_1^{*} \pi_{V}^{*} \mathcal F) = \Gamma(W, p_{1*} p_1^{*} \mathcal F_{G}).
$$
For the same reason, we also have an injection
$$
\Gamma(W, \varprojlim i^{-1} p_1^{(n)*} \mathcal F) \hookrightarrow \Gamma(W, \varprojlim p_1^{(n)*} \mathcal F_{G}).
$$
It is therefore sufficient to show that the morphism
$$
p_{1*} p_1^{*} \mathcal F_{G} \to \varprojlim p_1^{(n)*} \mathcal F_{G}
$$
is injective.
Now the question is local for the Grothendieck topology on $]X[_{P}$.
It is therefore also local on $P$ and we may assume that $P$ is affine, that $X$ is closed in $P$, and that we have local coordinates $t_{1}, \ldots, t_{n}$.
They induce local coordinates $\tau_{1}, \ldots, \tau_{n}$ on $P \times_{S} P$ with respect to the first projection $p_{1}$.
Using Lemma 4.4 of \cite {Berkovich99} (or proposition \ref{poly}), we get an isomorphism
$$
]X[_{P_{2}} \simeq ]X[_{\hat{\mathbf A}^n_{P}} = ]X[_{P} \times_{K} \mathbf B^n(0, 1^-).
$$
and we may use lemma \ref{rigproj} below. $\Box$\bigskip

\begin{lem} \label {rigproj}
Let $V$ be a good analytic variety and $\mathbf B$ any non trivial ball (open or closed) with coordinates $t_{1}, \ldots, t_{n}$.
Denote by
$$
p : V \times_{K} \mathbf B \to V
$$
and
$$
p^{(n)} : V \times_{K} \mathcal M( K[\underline t]/(\underline t)^n) \to V,
$$
the projections.
If $\mathcal F$ is a coherent $\mathcal O_{V}$-module, the canonical map
$$
p_{*}p^*\mathcal F \to \varprojlim_{n} p^{(n)*}\mathcal F
$$
is injective.
\end{lem}

\textbf{Proof:}
Since we work with coherent sheaves and good analytic spaces, we may freely use the Grothendieck topology.
Taking inverse limits, it is sufficient to consider the case of a closed ball. We may even assume that it has radius $1$.
Moreover, it is sufficient to check that, when $V$ is affine, the map
$$
\Gamma(V, p_{*}p^*\mathcal F) \to \Gamma(V, \varprojlim_{n} p^{(n)*}\mathcal F)
$$
is injective.
If we let $A := \Gamma(V, \mathcal O_{V})$ and $M := \Gamma(V, \mathcal F)$, this map is the canonical map
$$
M \otimes_{A} A\{\underline t\} \to M \otimes_{A} A[[\underline t]].
$$
By induction on the number of generators of $M$, we may assume that $M$ is a quotient of $A$.
It is therefore an affinoid algebra and we are reduced to the case $M = A$ which is part of the definition of the ring of convergent power series. $\Box$\bigskip

\begin{cor} \label{stratover} If $X \hookrightarrow P$ is a good admissible formal embedding into a quasi-compact formal $S$-scheme which is smooth at $X$, the canonical functor
$$
\textrm{Strat}^\dagger(X, P/S) \to \textrm{Strat}(X, P/S)
$$
is fully faithful on coherent modules.
\end{cor}

\textbf {Proof : }
Follows from the above description of morphisms in both categories. $\Box$\bigskip

 %
%
\section{Crystals and connections} \label{comp}

In this last section (as in the previous one), we fix a formal $\mathcal V$-scheme $S$.
\bigskip

If $(X, V)$ is an analytic variety over $S$ we will denote by $X_V/S$ the sieve generated by the (image of the) canonical morphism $(X,V) \to X/S$ of overconvergent presheaves on $\mathcal V$.
When $V = P_K$, we will write $X_P$ instead of $X_V$.
By definition, if $(X, V)$ is an analytic variety over $S$, the category $\textrm{AN}^\dagger(X_V/S)$ is essentially the full subcategory of $\textrm{AN}^\dagger(X/S)$ consisting of all $(X', V')$ such that the canonical map $X' \to X$ lifts to \emph{some} morphism $(X', V') \to (X, V)$.
Of course, any morphism $(f, u) : (X', V') \to (X, V)$ of analytic varieties over $S$ induces a morphism of analytic presheaves $f_u : X'_{V'}/S \to X_V/S$ which in turn gives a morphism of toposes
$$
f_{u} : (X_V/S)_{\textrm{AN}^\dagger} \to (X'_{V'}/S)_{\textrm{AN}^\dagger}.
$$
Since we will need it later, note also the following.
Let $(X, V)$ be an analytic variety over $S$ and $W$ an open subset of $S_K$. If $p_V : V \to S_K$ denotes the structural map and $V' := p_{V}^{-1}(W)$, then
$$
\Gamma(W, p_{X_V/S*}\mathcal F) = \Gamma(X_{V'}/S, \mathcal F).
$$
\begin{lem}
Let $f : X' \to X$ be a morphism of $S_k$-schemes that extends in two ways to morphisms
$$
(f, u_1) : (X', V') \to (X, V)
$$
and
$$
(f, u_2) : (X', V') \to (X, V)
$$
of analytic varieties over $S$.
If $E$ is an overconvergent crystal on $(X_V/S)$, there is a canonical isomorphism
$$
\epsilon_V : u_2^\dagger E_V \simeq u_1^\dagger E_V.
$$
\end{lem}

\textbf {Proof : }Just take
$$
\epsilon_V := \phi_{u_1}^\dagger \circ (\phi_{u_2}^\dagger)^{-1} : u_2^\dagger E_V \simeq E_{V'} \simeq u_1^\dagger E_V.
$$ $\Box$\bigskip

\begin{prop} \label{crystrat}
Let $(X, V)$ an analytic variety over $S$.
Then, the functor
$$
\textrm{Cris}^\dagger(X_V/S) \to \textrm{Strat}^\dagger(X, V/S), \quad E \mapsto (\mathcal F, \epsilon)
$$
with $\mathcal F := E_V$ and
$$
\epsilon :=\epsilon_V : p_2^\dagger E_V \simeq p_1^\dagger E_V
$$
is an equivalence of categories.
Moreover, we have
$$
p_{X_V/S*} E \simeq \mathcal H^{\dagger}(\mathcal F).
$$
\end{prop}

\textbf {Proof : }If $(f, u) : (X', V') \to (X, V)$ is a morphism in $\textrm{AN}(S)$, we have to set $E_V' := u^\dagger \mathcal F$. Now, if we are given two such morphisms $(f, u_1) : (X', V') \to (X, V)$ and $(f, u_2) : (X', V') \to (X, V)$, we can consider the diagonal morphism $(f, u) : (X', V') \to (X, V^2)$. Pulling back the Taylor isomorphism gives a canonical isomorphism $u_2^\dagger \mathcal F \simeq u_1^\dagger \mathcal F$. This shows that our definition is essentially independent of the choices. One easily checks that this gives a quasi-inverse to our functor.

We still have to prove the second assertion.
We already mentioned that if $W$ is an open subset of $S_K$ and $V' := p_{V/S_K}^{-1}(W)$, then
$$
\Gamma(W, p_{X,V/S*}E) = \Gamma(X_{V'}/S, E).
$$
Replacing $V$ by $V'$ , we are thus lead to show that $\Gamma(X_V/S, E)$ is the kernel of
$$
\epsilon \circ p_2^{-1} - p_1^{-1} : \Gamma(]X[_V, E_V) \to \Gamma(]X[_{V^{2}}, p_1^\dagger E_V)
$$
But we have
$$
\Gamma(X_V/S, E) = \varprojlim \Gamma((X', V'), E)
$$
$$
= \varprojlim \Gamma(]X'[_{V'}, E_{X',V'}) = \varprojlim \Gamma(]X'[_{V'}, ]f[_u^\dagger E_{X,V})
$$
Since any morphism factors through its graph, we see that
$$
\Gamma(X_V/S, E) = \varprojlim \left[
\Gamma(]X[_V, E_V)
\begin{array}{c}
\to \Gamma(]X[_{V^{2}}, p_1^\dagger E_V)\\
\to \Gamma(]X[_{V^{2}}, p_2^\dagger E_V)
\end{array}\right]
$$
which is what we want. $\Box$\bigskip

\begin{cor} \label {flatnesscr}
Let $(X, V)$ an analytic variety over $S$.
The realization functor to from $\textrm{Cris}^\dagger(X_V/S)$ to the category of $i_X^{-1}\mathcal O_V$-modules is right exact and faithful.
Moreover, if $V$ is universally flat in a neighborhood of $]X[_V$ over $S_K$, then $\textrm{Cris}^\dagger(X_V/S)$ is even an abelian category and the realization functor is exact.
\end{cor}
\textbf {Proof : } Follows from proposition \ref{flatness}. $\Box$
\bigskip

Note that there is a sequence of functors
$$\
\textrm{Cris}^\dagger(X_V/S) \simeq \textrm{Strat}^\dagger(X, V/S)
$$
$$
\to \textrm{Strat}(X, V/S) \to \textrm{MIC}(X, V/S).
$$
And there is also a sequence of morphisms
$$
p_{X_{V}/S*}E \simeq \mathcal H^\dagger(E_{V}) \to \varprojlim \mathcal H^{(n)}(E_{V}) \to E_V^{\nabla = 0}.
$$
In both cases, the last arrow is bijective when $\textrm{Char} K = 0$.

\begin{prop} \label{abel}
If $(X, V)$ is an analytic variety over $S$, then an overconvergent module $E$ on $X_V/S$ is finitely presented (resp. locally free of rank r) if and only if it is a crystal and $E_V$ is coherent (resp. locally free of rank r). Moreover, if $V$ is universally flat in the neighborhood of $]X[_{V}$ over $S_{K}$, then $\textrm{Mod}_{\textrm{fp}}^\dagger(X_V/S)$ is an abelian subcatgegory of $\textrm{Cris}^\dagger(X_V/S)$.
\end{prop}

\textbf {Proof : } This follows from corollary \ref{catT} because here, if $(f, u) : (X', V') \to (X, V)$ is any morphism of analytic varieties, then $E_{V'} = ]f[_{u}^\dagger E_{V}$ will also be coherent. $\Box$\bigskip

Note that if $E, E'$ are two overconvergent modules of finite presentation on $X_V/S$, then the dagger stratified module associated to $\mathcal H\textrm{om}_{\mathcal O_{X_{V}/S}^\dagger}(E,E')$ is $\mathcal H\textrm{om}_{i_{X}^{-1}\mathcal O_{V}}(E_{X,V}, E'_{X,V})$.

\begin{prop} Assume $\textrm{Char} K = 0$.
Let $X \hookrightarrow P$ be a good admissible formal embedding into a quasi-compact formal $S$-scheme which is smooth at $X$.
Then,
\begin{enumerate}
\item The functor $\textrm{Mod}^\dagger_{\textrm{fp}}(X_P/S) \to \textrm{MIC}(X, P/S)$ is fully faithful and its image consists of coherent modules with integrable connexion whose Taylor series converges on a strict neighborhood of $X$ in $P$.
\item If $E$ is an overconvergent module of finite presentation on $X_P/S$, we have
$$
p_{X_P/S*} E = E_{P}^{\nabla = 0}.
$$
\end{enumerate}
\end{prop}

\textbf {Proof : }
It follows from proposition \ref{crystrat} that $\textrm{Mod}^\dagger_{\textrm{fp}}(X_P/S)$ is equivalent to the category of coherent overconvergent stratified modules and we proved in corollary \ref{stratover} that the forgetful functor to stratified modules is fully faithful.
Finally, stratified modules are equivalent to modules with integrable connexions since we assumed $\textrm{Char} K = 0$. $\Box$\bigskip

\begin{cor} \label{rigthm}
Assume $\textrm{Char} K = 0$.
Let $X \hookrightarrow P$ be a good admissible formal embedding into a formal $S$-scheme which is smooth at $X$.
Then, we have a canonical equivalence of categories
$$
\textrm{Mod}^\dagger_{\textrm{fp}}(X_{P}/S) \simeq Isoc^\dagger(X \subset \bar X/S).
$$ 
\end{cor}

\textbf {Proof : } Recall that if $V$ is a paracompact analytic variety, and $V_{rig}$ denotes its set of rigid points with its rigid analytic topology, we have we have an equivalence of ringed toposes $\tilde V_{rig} \simeq \tilde V_{G}$. In particular, if $V$ is good, we get an equivalence of categories $\textrm{Coh}(\mathcal O_{V}) \simeq \textrm{Coh}(\mathcal O_{V_{rig}})$.
Therefore, it follows from propositions \ref{strictneigh} and \ref{jdagger} that $Isoc^\dagger(X \subset \bar X/S)$ is equivalent to the full subcategory of coherent $i^{-1} \mathcal O_{P_{K}}$-modules with an integrable connexion whose Taylor series converges on a strict neighborhood of the diagonal. $\Box$\bigskip

At last we will use the geometric results of section \ref {tech}.
Actually, we will only need Theorem \ref{main} that we can reformulate\ Êin the following way :

\begin{thm} \label{section}
Let $u : P' \to P$ be a morphism of good admissible formal embeddings over $S$ of an $S_k$-scheme $X$ which is proper and smooth at $X$. Then, the morphism
$$
(Id_X, u_K) : (X, P'_K) \to (X, P_K). 
$$
has locally a section in $\textrm{AN}^\dagger(S)$.
\end{thm} $\Box$ \bigskip

If $R$ is a sieve of an object $X$ in a site $C$, then $R$ is a covering sieve if and only if the induced morphism on associated sheaves $\tilde R \to \tilde X$ is an isomorphism. And in this case, we will obtain an equivalence $\tilde C_{/X} \simeq \tilde C_{/R}$.

\bigskip

\begin{cor} \label {final2}
With the assumptions and notations of the theorem, we have a canonical isomorphism of toposes
$$
(X_{P'}/S)_{\textrm{AN}^\dagger} \simeq (X_P/S)_{\textrm{AN}^\dagger}.
$$
\end{cor}

\textbf {Proof : }It follows from the theorem that the morphism of sieves of $X/S$, $X_{P'}/S \hookrightarrow X_P/S$, which is always a monomorphism, induces an isomorphism on the associated sheaves. Therefore, we get an equivalence on localized toposes. $\Box$\bigskip

\begin{prop} \label {final}			
Let $X \hookrightarrow P$ be a good admissible formal embedding with $P$ proper and smooth over $S$ at $X$, then $(X, P_K)$ is a covering of $X/S$ in $\textrm{AN}^\dagger(\mathcal V)$.
\end{prop}

\textbf {Proof : }We have to show that, given any $(U \subset Q \leftarrow V) \in \textrm{AN}^\dagger(X/S)$, there exists locally, a morphism to $(X \subset P \leftarrow P_K)$. It directly follows from the definition of the underlying category that it is sufficient to consider the case where $V = Q_K$. Moreover, the embedding $Q^{flat} \hookrightarrow Q$ of the maximal admissible formal subscheme gives an isomorphism on the generic fibers. We may therefore also assume that $Q$ is admissible.
We consider the graph $U \subset Q \times_S P$ and the projections
$$
p : Q \times_S P \to P, \quad q : Q \times_S P \to Q.
$$
It is sufficient to show that $(Id_U, q_K)$ has locally a section in $\textrm{AN}^\dagger(X/S)$ but this follows from the Theorem \ref {section} since by construction $q$ is proper and smooth at $U$. $\Box$\bigskip

\begin{cor} \label {restric}
If $X \hookrightarrow P$ is a good admissible formal embedding over $S$ with $P$ proper and smooth at $X$, we get an equivalence of toposes
$$
(X_P/S)_{\textrm{AN}^\dagger} \simeq (X/S)_{\textrm{AN}^\dagger}.
$$
\end{cor}

\textbf {Proof : } This follows from Proposition \ref {final} which tells us that $X_P$ is a covering sieve of $X$ and therefore, the sheaf associated to $X_P/S$ is exactly $X/S$. $\Box$\bigskip

\begin{prop} \label{retricris}
Let $X$ be an $S_k$-scheme.
\begin{enumerate}

\item If $u : P' \to P$ is a morphism of good admissible formal embeddings of $X$ over $S$ which is proper and smooth at $X$, we have a canonical equivalence
$$
\textrm{Cris}^\dagger(X_{P}/S) \simeq \textrm{Cris}^\dagger(X_{P'}/S).
$$
\item If $X \hookrightarrow P$ is a good admissible formal embedding over $S$ with $P$ proper and smooth at $X$, then we have an equivalence of categories
$$
\textrm{Cris}^\dagger(X/S) \simeq \textrm{Cris}^\dagger(X_P/S)
$$
\end{enumerate}

\end{prop}

\textbf {Proof : }The first assertion follows from Corollary \ref {final2} and the second one from Corollary \ref {restric}. $\Box$\bigskip

\begin{thm} Assume $\textrm{Char} K = 0$.
Let $X \hookrightarrow P$ be a good admissible embedding into a formal $S$-scheme which is proper and smooth at $X$.
Then,
\begin{enumerate}
\item The functor
$$
\textrm{Mod}^\dagger_{\textrm{fp}}(X/S) \to \textrm{MIC}(X, P/S), \quad E \mapsto (E_{P}, \nabla)
$$
is fully faithful and its image consists of modules with integrable connexion whose Taylor series converges on a strict neighborhood of $X$ in $P^2$.
\item If $E$ is an overconvergent module of finite presentation on $X/S$, we have
$$
p_{X/S*} E = E_{P}^{\nabla = 0}.
$$
\end{enumerate}
\end{thm}

\textbf {Proof : } This follows from theorem \ref{retricris}.
More precisely, we have the following sequence of functors
$$
\textrm{Cris}^\dagger(X/S) \simeq \textrm{Cris}^\dagger(X_P/S)
$$
$$
\simeq \textrm{Strat}^\dagger(X, P/S) \hookrightarrow \textrm{Strat}(X, P/S) \simeq \textrm{MIC}(X, P/S)
$$
and a sequence of morphisms
$$
p_{X/S*}E \simeq p_{X_{P}/S*}E \simeq \mathcal H^\dagger(E_{P}) \hookrightarrow \varprojlim \mathcal H^{(n)}(E_{P}) \simeq E_P^{\nabla = 0}.
$$
(in both cases, the last arrow also is bijective because $\textrm{Char} K = 0$). $\Box$\bigskip

For the last results, we consider only varieties which admit a good admissible formal embedding into a quasi-compact formal $S$-scheme which is proper and smooth at $X$.
This applies in particular to quasi-projective varieties.
In general, it would be necessary to glue, but as it is the tradition in rigid cohomology, we will not get through this matter here. We will consider the general case in a future article.

\begin{prop} If $X$ an $S_k$-scheme, then $\textrm{Mod}_{\textrm{fp}}^\dagger(X/S)$ is an abelian category.
Moreover, if $f : X' \to X$ is any $S_k$-morphism, then
$$
f_{AN}^{-1} : \textrm{Mod}_{\textrm{fp}}^\dagger(X/S) \to \textrm{Mod}_{\textrm{fp}}^\dagger(X'/S)
$$
is exact.
\end{prop}

\textbf {Proof : }
The first assertion results from proposition \ref{abel}.
For the second one, we may choose a good admissible embedding, both for $X$ and for $Y$, which is proper and smooth at $X$, respectively $Y$.
We may then use the diagonal embedding and assume that $f$ extends to a morphism of formal schemes which is smooth in the neighborhood of $X$.
In particular, it induces a flat morphism on strict neighborhoods, and therefore an exact functor. $\Box$\bigskip

\begin{prop} Assume $\textrm{Char} K = 0$.
If $X$ is an algebraic variety over $S_{k}$, we have a canonical equivalence of categories
$$
\textrm{Mod}^\dagger_{\textrm{fp}}(X/S) \simeq \textrm{Isoc}^\dagger(X/S)
$$
\end{prop}

\textbf {Proof : }
Follows from proposition \ref{retricris} and corollary \ref{rigthm}. $\Box$\bigskip

\textbf{Remark : } As a corollary, we recover the fact that the category $\textrm{Isoc}^\dagger(X/S)$ is essentially independant of the choices.

\bibliographystyle{plain}

\begin{thebibliography}{99}

\bibitem{Berkovich93}Vladimir G. Berkovich. \textsl{\'Etale cohomology for non-Archimedean analytic spaces}. Inst. Hautes \'Etudes Sci. Publ. Math., (78):5--161, 1993. 
\bibitem{Berkovich94} Vladimir G. Berkovich. \textsl{Vanishing cycles for formal schemes}. Invent. Math., 115(3):539--571, 1994. 
\bibitem{Berkovich99} Vladimir G. Berkovich. \textsl{Smooth p-adic analytic spaces are locally contractible}. Invent. Math., 137(1):1--84, 1999. 
\bibitem{Berthelot96*} Pierre Berthelot. \textsl{Cohomologie rigide et cohomologie \`a support propre, premi\`ere partie}. Pr\'epublication de lÕIRMAR, 96--03:1--89, 1996.
\bibitem{Berthelot96} Pierre Berthelot. \textsl{D-modules arithm\'etiques. I. Op\'erateurs diff\'erentiels de niveau fini}. Ann. Sci. \'Ecole Norm. Sup. (4), 29(2):185--272, 1996. 
\bibitem{Besser00} Amnon Besser. \textsl{A generalization of ColemanÕs p-adic integration theory}. Invent. Math., 142(2):397--434, 2000. 
\bibitem{Chiarellotto00} Bruno Chiarellotto. \textsl{Espaces de Berkovich et \'equations diff\'erentielles p-adiques. Une note}. Rend. Sem. Mat. Univ. Padova, 103:193--209, 2000. 
\bibitem{ChiarellottoTsuzuki03} Bruno Chiarellotto and Nobuo Tsuzuki. \textsl{Cohomological descent of rigid cohomology for \'etale coverings}. Rend. Sem. Mat. Univ. Padova, 109:63--215, 2003. 
\bibitem{Ducros03*} Antoine Ducros. \textsl{\'Etude de certaines propri\'et\'es locales et globales des espaces de Berkovich}. Pr\'epublication de lÕIRMAR, 03--41:1Ð59, 2003. 
\bibitem{EtesseLeStum93} Jean-Yves \'Etesse and Bernard Le Stum. \textsl{Fonctions L associ\'ees aux F -isocristaux surconvergents I. Interpr\'etation cohomologique}. Math. Ann., 296(3):557--576, 1993. 
\bibitem{GabrielZisman67} P. Gabriel and M. Zisman. Calculus of fractions and homotopy theory. Ergebnisse der Mathematik und ihrer Grenzgebiete, Band 35. Springer-Verlag New York, Inc., New York, 1967. 
\bibitem{KashiwaraSchapira90} Masaki Kashiwara and Pierre Schapira. Sheaves on manifolds, volume 292 of Grund lehren der Mathematischen Wissenschaften. Springer-Verlag, Berlin, 1990.
\bibitem{Temkin00} Michael Temkin. \textsl{On local properties of non-Archimedean analytic spaces}. Math. Ann., 318(3):585--607, 2000.
\bibitem{Temkin04} Michael Temkin. \textsl{On local properties of non-Archimedean analytic spaces. II}. Israel J. Math., 140:1--27, 2004.
\end{thebibliography}

\Addresses
\end{document}